\documentclass[reqno,a4paper]{amsart}

\usepackage{mathrsfs,amsmath,amssymb,graphicx}
\usepackage{stmaryrd}
\usepackage[unicode]{hyperref}
\usepackage{latexsym}
\usepackage{layout}
\usepackage[english]{babel}
\usepackage{color}
\usepackage{fancyvrb}
\usepackage{subfigure}

\theoremstyle{plain}
\newtheorem{theorem}{Theorem}[section]

\newtheorem{proposition}[theorem]{Proposition}

\newtheorem{definition}{Definition}[section]
\theoremstyle{remark}
\newtheorem{remark}{Remark}[section]

\newcommand{\Z} {\mathbb Z}
\newcommand{\SSS} {\mathcal S}
\newcommand{\sla} {\boxslash}
\newcommand{\bsla} {\boxbslash}

\graphicspath{{figures/}}

\title{Exploring the ``Rubik's Magic'' universe}

\author{Maurizio Paolini}
\address{
Dipartimento di Matematica e Fisica,
Universit\`a Cattolica ``Sacro Cuore'',
Brescia, Italy
E-mail: paolini@dmf.unicatt.it}

\date{\today}

\begin{document}


\begin{abstract}
By using two different invariants for the Rubik's Magic puzzle, one of metric
type, the other of topological type, we can dramatically reduce the universe
of constructible configurations of the puzzle.
Finding the set of actually constructible shapes remains however a challenging task,
that we tackle by first reducing the target shapes to specific configurations:
the octominoid 3D shapes, with all tiles parallel to one coordinate plane; and the
planar ``face-up'' shapes, with all tiles (considered of infinitesimal width) lying in 
a common plane and without superposed consecutive tiles.
There are still plenty of interesting configurations that do not belong to either
of these two collections.
The set of constructible configurations (those that can be obtained by manipulation of
the undecorated puzzle from the starting situation) is a subset of the set of configurations
with vanishing invariants.
We were able to actually construct all octominoid shapes with vanishing invariants
and most of the planar ``face-up'' configurations.
Particularly important is the topological invariant, of which we recently found mention in
\cite{Ver:87} by Tom Verhoeff.
\end{abstract}

\maketitle

\section{Introduction}\label{sec:intro}
The \textit{Rubik's Magic} is another creation of Ern\H{o} Rubik, the brilliant hungarian inventor
of the ubiquitous ``cube'' that is named after him.
The \textit{Rubik's Magic} puzzle is much less known and not very widespread today,
however it is a really surprising object that hides aspects which makes it quite an
interesting subject for a mathematical analysis on more than one level.

We investigate here two different invariants that can
be used to prove the unreachability of many spatial configurations of the puzzle, one of
these invariants, of topological type, is to our knowledge never been extensively studied before,
although it is presumably the same mentioned in \cite{Ver:87}, and allows to
significantly reduce the number of theoretically constructible shapes.
However even in the special case of planar ``face-up'' configurations (see
Section \ref{sec:piane}) we don't know whether the combination
of the two invariants, together with basic constraints coming from the mechanics of the
puzzle, is complete, {\it i.e.} if it characterizes the set of constructible configurations.
Indeed there are still a few planar face-up configurations having both vanishing invariants,
but that we are not able to construct.
In this sense this Rubik's invention remains an interesting subject of mathematical analysis.

In Section \ref{sec:rompicapo} we describe the puzzle and discuss its mechanics, the local
constraints are discussed in Section \ref{sec:vincolilocali}.
The addition of a ribbon (Section \ref{sec:nastro}) allows to introduce the two invariants,
the metric and the topological invariants, described respectively in Sections
\ref{sec:invariantemetrico} and \ref{sec:invariantetopologico}.

The set of \textit{octominoid} shapes (all tiles are parallel to a coordinate plane and no two of them
are superposed) is described in Section \ref{sec:octominoid}.
There are a total of $460$ distinct
octominoid configurations of the undecorated puzzle with vanishing invariants and \textbf{all} of them are actually
contructible with the real puzzle \cite[3D octominoids]{mine:web} meaning that within this special class of shapes
the two invariants are complete.

The special ``face-up'' planar configurations are defined in Section \ref{sec:piane} and
their invariants computed in Section \ref{sec:forsecostruibili}.
There are a total of $25$ configurations with vanishing invariants, all of which we were
able to actually construct \cite[planar face-up configurations]{mine:web} with only $5$ exceptions
(Section \ref{sec:nonclassificate}).
The two basic configurations of Figures \ref{fig:black2x4} and \ref{fig:blacksolved} are contained
in both octominoid and planar face-up classes, for a total of $485$ configurations.
Of course there are still many configurations that are not contained in either of the two
classes (actually there are infinitely many of them), making the exploration of the Rubik's
Magic universe far from complete.

%
%

We conclude the paper with a brief description of the software codes used to help
in the analysis of the octominoid and of the planar face-up configurations (Section \ref{sec:code}).

\section{The puzzle}\label{sec:rompicapo}

\begin{figure}\begin{center}
\includegraphics[width=7cm]{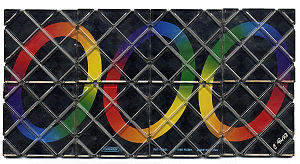}
\includegraphics[width=5cm]{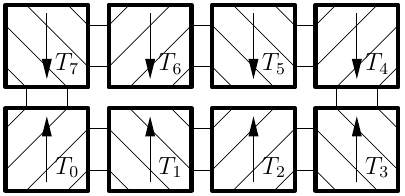}
\caption{\small{The original puzzle in its starting configuration (left).
Orientation scheme for the tiles (right).}}
\label{fig:black2x4}
\end{center}\end{figure}

\begin{figure}\begin{center}
\includegraphics[width=6cm]{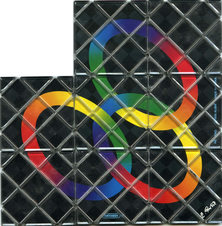}
\caption{\small{The puzzle in its target configuration, the tiles
are turned over with respect to Figure
\ref{fig:black2x4}.}}
\label{fig:blacksolved}
\end{center}\end{figure}

The \textit{Rubik's Magic} puzzle (see Figure \ref{fig:black2x4} left)
consists of $8$ decorated square tiles positioned to
form a $2\times 4$ rectangle.

They are ingeniously connected to each other by means of nylon strings
lying in grooves carved on the tiles and tilted $45$ degrees
\cite{basteleien:web}.

The tiles are decorated in such a way that on one side of the $2\times 4$
original configuration we can see the picture of three unconnected rings,
whereas on the back side we find non-matching drawings representing parts
of rings with crossings among them.

The declared aim is to manipulate the tiles in order to correctly place
the decorations on the back, which can be done only by changing the global
outline of the eight tiles.

The solved puzzle is shown in Figure \ref{fig:blacksolved} with the tiles
positioned in a $3 \times 3$ square with a missing corner and overturned
with respect to the original configuration of Figure \ref{fig:black2x4}.

Detailed instructions on how the puzzle can be solved and more generally
on how to construct interesting shapes can be copiously found in the
internet, we just point to the Wikipedia entry \cite{wikipedia:web} and to
the web page \cite{basteleien:web}.
The booklet \cite{Nou:86} contains a detailed description of the puzzle and
illustrated instructions on how to obtain particular configurations.

The decorations can be used to distinguish a ``front'' and a ``back'' face
of each tile and to orient them by suitably chosing an ``up'' direction.

After dealing with the puzzle for some time it becomes apparent that a
few local constraints are always satisfied.
In particular the eight tiles remain always connected two by two in such
a way to form a cyclic sequence.
To fix ideas let us denote the eight tiles by $T_0$, $T_1$, ..., $T_7$,
with $T_0$ the lower-left tile in Figure \ref{fig:black2x4} and the others
numbered in counterclockwise order.
For example tile $T_3$ is the one with the Rubik's signature in its
lower-right corner (see Figure \ref{fig:black2x4}).

With this numbering tile $T_i$ is always connected through one of its sides
to both tiles $T_{i+1}$ and $T_{i-1}$.
Here and in the rest of this paper we shall always assume the index
$i$ in $T_i$ to be defined ``modulo $8$'', {\it i.e.} that for example
$T_8$ is the same as $T_0$.

We shall conventionally orient the tiles such that in the initial configuration
of Figure \ref{fig:black2x4} all tiles are ``face up'' (i.e. with their
front face visible), the $4$ lower tiles ($T_0$ to $T_3$) are ``straight''
(not upside down), the $4$ upper tiles ($T_4$ to $T_7$) are ``upside down''
(as a map with the north direction pointing down),
see Figure \ref{fig:black2x4} right.

At a more accurate examination it turns out that only half of the grooves
are actually used (those having the nylon threads in them).
These allow us to attach to a correctly oriented tile (face up and straight)
a priviledged direction:
direction $\sla$ (``slash'': North-East to South-West) and direction
$\bsla$ (``backslash'': North-West to South-East).
The used groves are shown in Figure \ref{fig:black2x4} right.
From now on we shall disregard completely the unused grooves.
In the initial configuration tiles $T_i$ with even $i$ are all tiles
of type $\sla$, whereas if $i$ is odd we have a tile of type $\bsla$.

The direction of the used grooves in the back of a tile is opposite
(read orthogonal) to the direction of the used grooves of the front
face,
but beware that when we revert (turn over) a tile a $\sla$ groove becomes
$\bsla$, so that the reversed tile remains of the same type
($\sla$ or $\bsla$).

From the point of view of an idealized physical modelling a natural choice would be to
assume that the
tiles are made of a rigid material and have infinitesimal thickness,
and that the nylon threads are perfectly flexible but inextensible (and
of infinitesimal thickness).
This allows for two or more tiles to be juxtaposed in space, however in such
a case we still need to retain the information about their relative position
(which is above which).

However in this model there are moves that can be performed on the real puzzle (that entail a small
amount of elongation on the nylon wires) but that are \textbf{not} allowed
in the ideal model (see e.g. the two interesting shapes denoted \textit{armchair}
and \textit{hard-to-reach planar shape} linked from \cite{mine:web}.

For this reason in the real constructions we shall allow for moves that
entail a small amount of deformation of the tiles and elongation of the wires.
We shall not be rigorous about what is allowed and what is not, the
rule being that we shall generally allow for moves that can actually be performed
on the real puzzle.


On the contrary the real puzzle has
non-infinitesimal tile thickness, which can lead to configurations that
are allright for the idealized physical model but that are difficult or impossible
to achieve (because of the imposed stress on the nylon threads) with
the real puzzle.


\subsection{Undecorated puzzle}

We are here mainly interested in the study of the \emph{shapes} in space
that can be obtained,
so we shall neglect the decorations on the tiles and only consider the
direction of the grooves containing the nylon threads.
In other words we only mark one diagonal on each face
of the tiles, one connecting two opposite vertices on the front face and
the other connecting the remaining two vertices on the back face.

Now the tiles (marked with these two diagonals) are indistinguishable;
distinction between $\sla$ and $\bsla$ is only possible after we have
``oriented'' a tile and in such a case rotation of $90$ degrees or a mirror reflection
will exchange $\sla$ with $\bsla$.

\begin{definition}[orientation]\label{def:orientation}
A tile can be oriented by drawing on \textbf{one} of the two faces an
arrow parallel to a side.
We have thus eight different possible orientations.
We say that two adjacent tiles are compatibly oriented  if their arrows
perfectly fit together (parallel and pointing to the same direction) when we
ideally rotate one tile around the side on which they are hinged to make it
juxtaposed to the other.
There is exactly one possible orientation of a tile that is compatible with the
orientation of an adjacent tile.
A configuration of tiles is \textbf{orientable} if it is possible to orient
all tiles such that they are pairwise compatibly oriented.
For an orientable configuration we have eight different choices for a compatible
orientation of the tiles.
\end{definition}

An example of compatible orientation of a configuration is shown in Figure
\ref{fig:black2x4} right, which makes the initial $2 \times 4$ configuration 
orientable.
Once we have a compatible orientation for a configuration, we can classify
each tile as $\sla$ or $\bsla$ according to the relation between the orienting
arrow and the marked diagonal:
a tile is of type $\sla$ if the arrow alignes with the diagonal after a clockwise
rotation of $45$ degrees,
it is of type $\bsla$ if the arrow alignes with the diagonal after a counterclockwise
rotation of $45$ degrees.
Two adjacent tiles are always of opposite type.

\begin{definition}\label{def:chiral}
A spatial configuration of the puzzle that is \textbf{not}
congruent (also considering the marked diagonals) after a rigid motion
with its mirror image will be called \textbf{chiral},
otherwise it will be called \textbf{achiral}.
Note that a configuration is achiral if and only if it is mirror
symmetric with respect to some plane.
\end{definition}

The initial $2 \times 4$ configuration is achiral since it is specularly
symmetric with respect to a plane orthogonal to the tiles.

\begin{definition}
We say that an orientable spatial configuration of the puzzle (without
decorations) \textbf{is constructible} if it can be obtained from the
initial $2 \times 4$ configuration through a sequence of admissible
moves of the puzzle.
\end{definition}

Once we have identified all the constructible spatial configurations,
we also have all constructible configurations of the decorated puzzle.
This is a consequence of the fact that all possible $2 \times 4$
configurations of the decorated puzzle are completely classified
(see for example \cite{basteleien:web} or \cite{Nou:86}).

We note here that all $2 \times 4$ configurations of the undecorated
puzzle are congruent, however the presence of the marked diagonal might
require a reversal of the whole configuration upside-down
in order to obtain the congruence.

For chiral configurations (those that cannot be superimposed with their
specular images) the following result is useful.

\begin{theorem}\label{teo:mirrorconstructible}
A spatial configuration of the undecorated puzzle is constructible if and
only if its mirror image is constructible.
\end{theorem}

\begin{proof}
If a configuration is constructible we can reach it by a sequence of moves
of the puzzle starting from the initial $2 \times 4$ configuration.
However the initial $2 \times 4$ configuration is specularly symmetric,
hence we can perform the specular version of that sequence of moves to
reach the specular image of the configuration that we are considering.
\end{proof}


\section{Local constraints}\label{sec:vincolilocali}
We now consider a version of the puzzle where in place of the usual decoration
we draw arrows on the ``front'' face of the tiles as in Figure
\ref{fig:black2x4} right.
The linking mechanism with the nylon threads is such that two consecutive
tiles $T_i$ and $T_{i+1}$ are always ``hinged'' together through one of their
sides.
In particular, if we suitably orient $T_i$ with its ``front'' face visible
and ``straight'', i.e. with the arrow visible and pointing up) and we rotate tile
$T_{i+1}$ such that its center is as far away as possible from the
center of $T_i$ (like an open book), then also $T_{i+1}$ will have its
arrow visible and
\begin{itemize}
\item
 \textbf{pointing up} if the two tiles are hinged through a vertical
side (the right or left side of $T_i$);
\item
 \textbf{pointing down} it the two tiles are hinged through a horizontal
side (the top or bottom side of $T_i$).
\end{itemize}

The surprising aspect of the puzzle is that when we ``close the book'',
i.e. we rotate $T_{i+1}$ so that it becomes superimposed with
(stacked on or below) $T_i$,
we than can ``reopen the book'' with respect to a different hinging
side.
The new hinging side is one of the two sides that are orthogonal to
the original hinging side, which one depending on the type of the
involved tiles (direction of the marked diagonals) and can be identified
by the rule that the new side is not separated from the previous one
by the ``inner'' marked diagonals.
For example, if $T_i$ is of type $\sla$ (hence $T_{i+1}$ is of type
$\bsla$) and they are hinged through the right side of $T_i$
(as $T_0$ and $T_1$ of Figure \ref{fig:black2x4} right) then after
closing the tiles by rotating $T_{i+1}$ \textbf{up} around its left side and
placing it on top of (stacked above) $T_i$, then we can reopen
the tiles with respect to the bottom side.
On the contrary, if we rotate \textbf{down} $T_{i+1}$, so that it
becomes stacked below $T_i$ (and the involved marked diagonal of $T_i$
is the one on the back face), the new hinging side will be the upper
side.

We remark that if a configuration does not contain superimposed consecutive
tiles, then the hinging side of any pair of consecutive tiles is uniquely
determined.
If the tiles are (compatibly) oriented, than for each tile $T_i$ we have a
unique side (say East, North, West or South, in short $E$, $N$, $W$ or $S$)
about which it is hinged with the preceding tile $T_{i-1}$ and a
unique side ($E$, $N$, $W$ or $S$) about which it is hinged with $T_{i+1}$.
The two sides can be the same.

\begin{definition}\label{def:localshape}
For a given spatial oriented configuration of the (undecorated) puzzle without
stacked consecutive tiles we say that a tile is
\begin{description}
\item[straight] if the two hinging sides are opposite;
\item[curving] if the two hinging sides are adjacent (but not the same).
In this case we can distinguish between tiles \textbf{curving left} and
tiles \textbf{curving right} with the obvious meaning and taking into
account the natural ordering of the tiles induced by the tile index;
\item[a flap] if it is hinged about the same side with both the previous
and the following tile.
\end{description}
\end{definition}

\subsection{Flaps}\label{sec:flaps}

Flap tiles (those that, following Definition \ref{def:localshape},
have a single hinging side with the two adjacent tiles) require a
specific analysis.
The term ``flap'' is the same used in \cite{Nou:86} and refers to the
similarity with the flaps of an airplane, that can rotate about a single
side.

Given an oriented configuration with a flap $T_i$, let us fix the attention to
the three consecutive tiles $T_{i-1}$, $T_i$, $T_{i+1}$ and ignore all the
others.
Place the configuration so that $T_i$ is horizontal, with its front face
up and the arrow pointing North, then rotate $T_{i-1}$ and $T_{i+1}$ at
maximum distance from $T_i$ so that they become reciprocally superimposed.

Now all three tiles have their front face up and we can distinguish
between two situations:

\begin{definition}\label{def:flaps}
Tile $T_i$ is an \textbf{ascending} flap if tile $T_{i-1}$ is \textbf{below}
tile $T_{i+1}$;
it is a \textbf{descending} flap in the opposite case.
Tile $T_i$ is a horizontal ascending/descending flap if it is hinged at a vertical
side (a side parallel to the arrow indicating the local orientation of the flap tile),
it is a vertical ascending/descending flap otherwise.
\end{definition}

\section{The ribbon trick}\label{sec:nastro}
In order to introduce the metric and the topological invariants we resort to a simple
expedient:
we insert a ribbon in between the tiles that more or less follows the path
of the nylon threads.

The ribbon is colored red on one side (front side) and blue in its back side
and is oriented with longitudinal arrows printed along its length that allows
to follow it in the positive or negative direction.

\begin{figure}
\centering
\includegraphics[height=4cm]{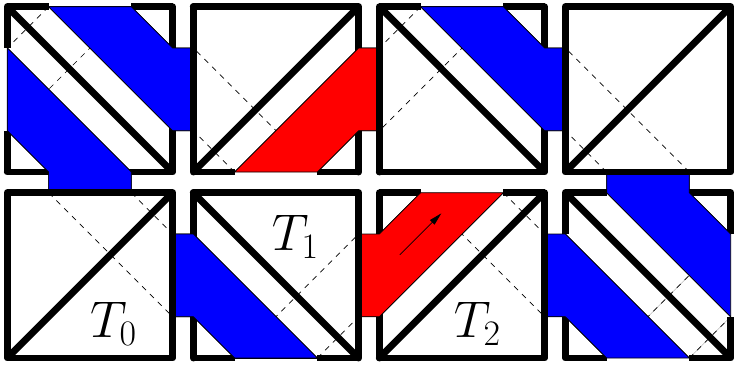}
\caption{\small{Ribbon path among the tiles.}}
\label{fig:nastro}
\end{figure}

Let the tiles have side of length $1$, then the ribbon has width that does not
exceed $\frac{\sqrt{2}}{4}$ (the distance between two nearby grooves), so
that it will not interfere with the nylon threads.
We insert the ribbon as shown in Figure \ref{fig:nastro}.
More precisely take the $2 \times 4$ initial configuration of the puzzle and
start with tile $T_2$.
Position the ribbon such that it travels diagonally along the front face of
$T_2$ as shown in Figure \ref{fig:nastro}, then wrap the ribbon around
the top side of $T_2$ and travel downwards along the back of $T_2$ to
reach the right side.
At this point we move from the back of $T_2$ to the front of $T_3$
(the ribbon now has its blue face up) and continue downward until we reach
the bottom side of $T_3$, wrap the ribbon on the back and so on.

In general, every time that the ribbon reaches a side of a tile that is not
a hinge side with the following tile, we wrap it around the tile (from the
front face to the back face or from the back face to the front face) as if
it ``bounces'' against the side.
Every time the ribbon reaches a hinging side of a tile with the following
tile it moves to the next tile and crosses from the back (respectively front)
side of one tile to the front (respectively back) side of the other and
maintains its direction.

In all cases the ribbon travels with sections of length $\delta = \frac{\sqrt{2}}{2}$
between two consecutive ``touchings'' of a side.
It can stay adjacent to a given tile during one, two or three of such $\delta$
steps:
one or three if the tile is a \textit{curving} tile (Definition \ref{def:localshape}), two if the tile is a \textit{straight} tile.

After having positioned the ribbon along all tiles, it will close on itself
nicely (in a straight way and with the same orientation) on the starting
tile $T_2$, and we tape it with itself.
In this way the total length of the ribbon is $16 \delta$ with an average of
$2\delta$ per tile,
moreover if we remove the ribbon without cutting it (by making the tiles
``disappear''), we discover that we can deform it in space into the lateral
surface of a large and shallow cylinder with height equal to the ribbon
thickness and circumference $16 \delta$.

Direct inspection also shows that the inserted ribbon does not impact on the
possible puzzle moves, whereas its presence allows us to define the two
invariants of Sections \ref{sec:invariantemetrico} and \ref{sec:invariantetopologico}.

We remark a few facts:

\begin{enumerate}
\item
The ribbon is oriented: it has arrows on it pointing in the direction in which
we have inserted it, and while traversing the puzzle along the ribbon the tiles
are encountered in the order given by their index.
\item
Each time the ribbon ``bounces'' at the side of a tile (moving from the
front face to the back face or viceversa) its direction changes of $90$ degrees
and simultaneously it turns over.
This does not happen when the ribbon moves from one tile to the next, it does not
change direction and it does not turn over.
\item
Each $\delta$ section of the ribbon connects a horizontal side to a vertical
side or viceversa;
consequently the ribbon touches alternatively horizontal sides and vertical
sides.
\item
Each time the ribbon touches a lateral side it goes from one side of the tiles
to the other (from the front to the back or from the back to the front).
\end{enumerate}

The above points 3 and 4 prove the following

\begin{proposition}
Following the orientation of the ribbon, when the ribbon touches\slash{}crosses a
vertical side, it ``emerges'' from the back of the tiles to the front,
whereas when it touches/crosses a horizontal side, it ``submerges'' from the
front to the back.
Here vertical or horizontal refers to the local orientation assigned to the tiles.
\end{proposition}


\subsection{Behaviour of the ribbon at a flap tile.}\label{sec:nastroflap}

It is not obvious how the ribbon behaves at a flap tile (such tiles are not present in the
initial $2 \times 4$ configuration).
We can reconstruct the ribbon position by imagining a movement that transforms a configuration
without flaps to another with one flap.

\begin{figure}\begin{center}
\includegraphics[height=2cm]{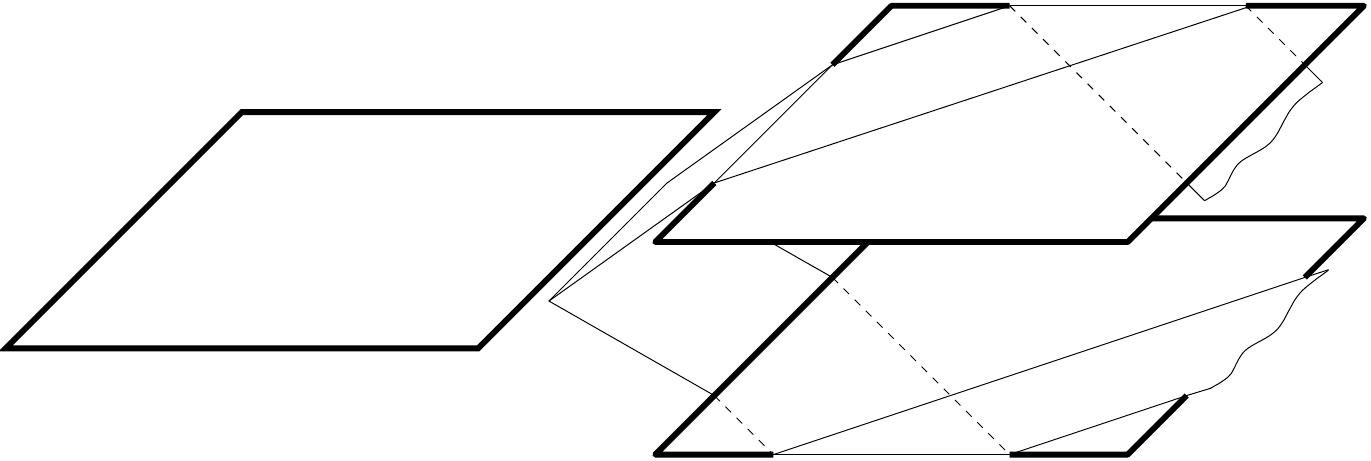}
\includegraphics[height=2cm]{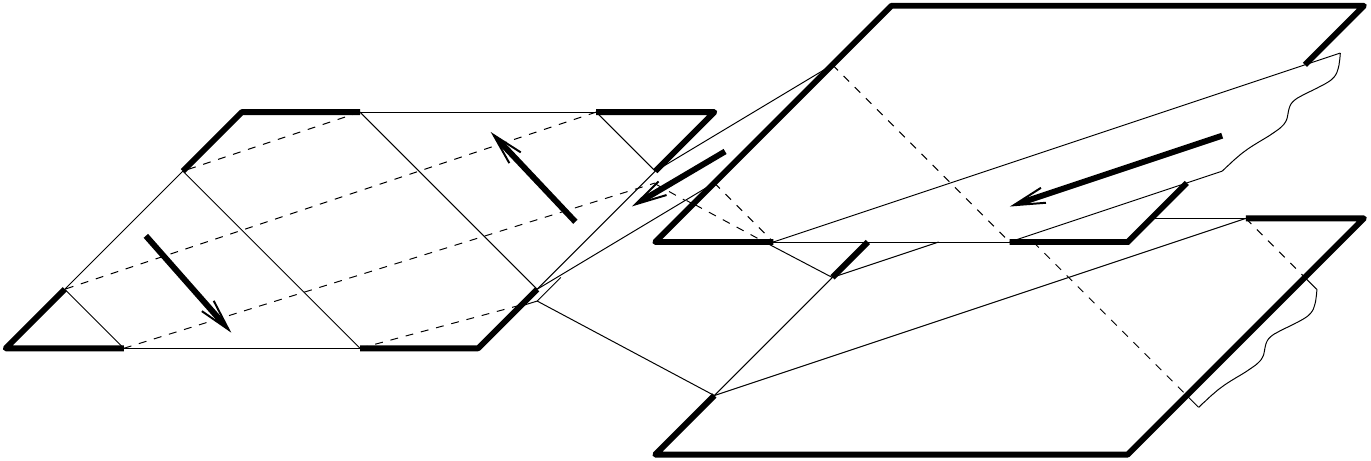}
\caption{\small{Position of the ribbon in presence of a ``flap'' tile.
The flap tiles are of type $\bsla$, The superposed tiles are all of type $\sla$.
Left: ascending flap, the ribbon does not even touch the flap tile.
Right: descending flap, the ribbon completely wraps the flap tile with
four sections, two on the upper (front) face and two on the back face.}}
\label{fig:nastroflap}
\end{center}\end{figure}

It turns out that there are two different situations.
In one case the ribbon completely
avoids to touch the flap tile $T_i$ and directly goes from $T_{i-1}$ to $T_{i+1}$,
this happens when in a neighbourhood of the side where the flap tile is hinged the
ribbon is on the front face of the upper tile and on the bottom face of the lower
tile (in the configuration where $T_{i-1}$ and $T_{i+1}$ are furthest away from $T_i$,
hence superposed),
this situation is illustrated in Figure \ref{fig:nastroflap} left.
In the other case the ribbon wraps around $T_i$ with four $\delta$ sections alternating
between the front face and the back face,
this situation is illustrated in Figure \ref{fig:nastroflap} right.

The first of the two cases arises at an \textit{ascending} flap hinged at a vertical
side (horizontal ascending flap) or at a \textit{descending} flap hinged at a horizontal
side (vertical descending flap);
this is independent of the type $\sla$ or $\bsla$ of the flap tile.

The second of the two cases arises at a vertical ascending flap or at a horizontal descending
flap.

\section{Metric invariant}\label{sec:invariantemetrico}
Whatever we do to the puzzle (with the ribbon inserted) there is no way to change the
length of the ribbon!

This allows to regard the length of the ribbon associated to a given spatial configuration
as an invariant, it cannot change under puzzle moves.
The computation of the ribbon length can be carried out by following a few simple
rules, they can also be found in \cite{Nou:86}.

The best way to proceed is to compute for each tile $T_i$ how many $\delta$ sequences of the ribbon
wrap it and subtract the mean value $2$.
The resulting quantity will be called $\Delta_i$ and its value is:
\begin{itemize}
\item $\Delta_i = 0$ if $T_i$ is a straight tile (Definition \ref{def:localshape});
\item $\Delta_i = -1$ if $T_i$ is of type $\sla$ and is ``curving left'', or if it is of
type $\bsla$ and is ``curving right'';
\item $\Delta_i = +1$ if $T_i$ is of type $\sla$ and is curving right, or if it is of type $\bsla$
and is curving left;
\item $\Delta_i = -2$ if $T_i$ is a horizontal ascending flap (Definition \ref{def:flaps}) or
a vertical descending flap (see Figure \ref{fig:nastroflap} left);
\item $\Delta_i = +2$ if $T_i$ is a horizontal descending flap or a vertical ascending flap
(see Figure \ref{fig:nastroflap} right).
\end{itemize}

The last two cases ($|\Delta_i| = 2$) follow from the discussion in Section \ref{sec:nastroflap}.

We call $\Delta = \sum_{i=0}^7 \Delta_i$, the sum of all these quantities, then the
total length of the ribbon will be $16 \delta + \Delta \delta$ and hence $\Delta$ is invariant
under allowed movements of the puzzle.
Since in the initial configuration we would have $\Delta = 0$ it follows that

\begin{theorem}\label{teo:invariantemetrico}
Any constructible configuration of the puzzle necessarily satisfies
$\Delta = 0$.
\end{theorem}

This invariant can also be found in \cite[page 19]{Nou:86}, though it is not actually
justified.

A few configurations (e.g. the $3 \times 3$ shape without the central square, called
``window shape'' in \cite{Nou:86}) can be ruled out as non-constructible by computing the
$\Delta$ invariant.
The ``window shape'' has a value $\Delta = \pm 4$, the sign depending on how we orient
the tiles.
It is non-constructible because $\Delta \neq 0$.
Figure \ref{fig:nc_shapes} (left) shows a deformed version of this shape.

Another interesting configuration that can be ruled out using this invariant is
sequence \eqref{eq:seqbirillo}, to be discussed in Section \ref{sec:noflaps}.


\section{Topological invariant}\label{sec:invariantetopologico}
Sticking to the ribbon idea (Section \ref{sec:nastro}) we seek a way to know
whether a given ribbon configuration (with the tiles and nylon threads removed)
can be obtained by deformations in space starting from the configuration where
the ribbon is the lateral surface of a cylinder.

Topologically the ribbon is a surface with a boundary, its boundary consists
of two closed strings.

One thing that we may consider is the center line of the ribbon: it is a
single closed string that can be continuously deformed in space and is not
allowed to cross itself.
Mathematically we call this a ``knot'', a whole branch of Mathematics is
dedicated to the study of knots, one of the tasks being finding ways to identify
``unknots'', i.e. tangled closed strings that can be ``unknotted'' to a perfect
circle.

This is precisely our situation: the center line of the ribbon must be an
unknot, otherwise the corresponding configuration of the puzzle cannot be
constructed.
However we are not aware of puzzle configurations that can be excluded for this
reason.

Another (and more useful) idea consists in considering the two strings forming the
boundary of the ribbon.
In Mathematics, a configuration consisting in possibly more than one closed string
is called a ``link''.
Here we have a two-components link that in the starting configuration can be
deformed into two unlinked perfect circles.

There is a topological invariant that can be easily computed, the \textit{linking
number} between two closed strings, that does not change under continuous deformations
of the link (again prohibiting selfintersections of the two strings or intersections
of one string with the other).

In the original configuration of the puzzle, the two strings bordering the ribbon
have linking number zero: it then must be zero for any constructible configuration.


\subsection{Computing the linking number}

\begin{figure}\begin{center}
\includegraphics[height=2cm]{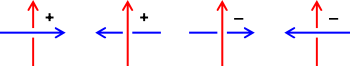}
\caption{\small{Signature of a crossing for the computation of the linking number.}}
\label{fig:linkingnumber}
\end{center}\end{figure}

In the field of \textit{knot theory} a knot, or more generally a link, is often
represented by its \text{diagram}.
It consists of a drawing on a plane corresponding to some orthogonal projection
of the link taken such that the only possible selfintersections are transversal
crossings where two distinct points of the link project onto the same point.
We can always obtain such a \textit{generic} projection possibly by changing
a little bit the projection direction.
We also need to add at all crossings the information of what strand of the
link passes ``over'' the other.
This is usually done by inserting a small gap in the drawing of the strand
that goes below the other, see Figure \ref{fig:linkingnumber}.

\begin{figure}\begin{center}
\includegraphics[height=3cm]{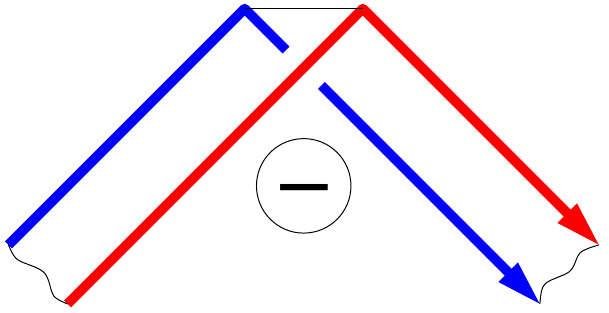}~~
\includegraphics[height=3cm]{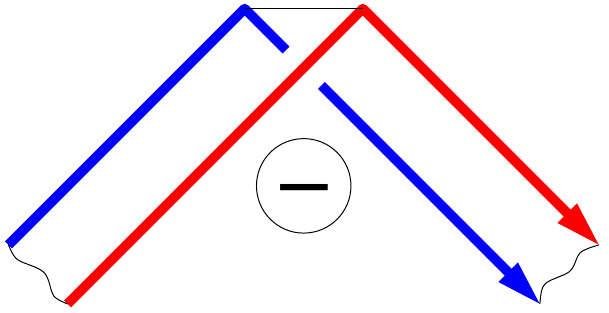}
\caption{\small{The ribbon \textit{bounces} at the side of a tile.}}
\label{fig:twistsign}
\end{center}\end{figure}

\begin{figure}\begin{center}
\includegraphics[height=3cm]{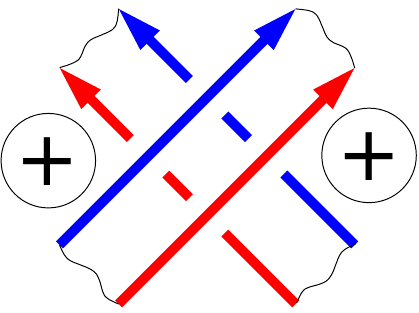}~~~~~~
\includegraphics[height=3cm]{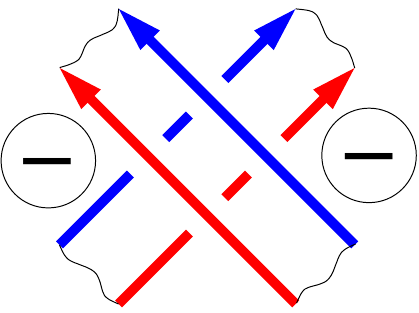}
\caption{\small{The ribbon passes over/below itself.}}
\label{fig:crosssign}
\end{center}\end{figure}

In order to define the linking number between two closed curves we need to select
an orientation (a traveling direction) for the two curves.
In our case the orientation of the ribbon induces an orientation of the two
border strings by following the same direction.
The linking number changes sign if we revert the orientation of one of the two curves,
so that it becomes insensitive upon the choice of orientation of the ribbon.
Once we have an orientation of the two curves, we can associate a signature
to each crossing as shown in Figure \ref{fig:linkingnumber} and a corresponding
weight of value $\pm \frac{1}{2}$.
Crossings of a component with itself are ignored in this computation.

The linking number is given by the sum of all these contributions.
Since the number of crossings in between the two curves in the diagram is necessarily
even, it follows that the linking number is an integer and it can be proved that it
does not change under continuous deformations of the link in space.
Two far away rings have linking number zero, two linked rings have linking number
$\pm 1$.

In our case we shall investigate specifically the case where all tiles are horizontal
and ``face-up'', in which case we have two different situations that produce
crossings between the two boundary strings.
We shall then write the linking number as the sum of a ``twist'' part ($L_t$) and
a ``ribbon crossing'' part ($L_c$)
\begin{equation}
L = L_t + L_c
\end{equation}
where we distinguish the two cases:

\begin{enumerate}
\item
The ribbon wraps around one side of a tile (Figure \ref{fig:twistsign}).
This entail one crossing in the diagram, that we shall call ``twist crossing'' since
it is actually produced by a twist of the ribbon.
A curving tile (as of Definition \ref{def:localshape}) can contain only zero or
two of this type of crossings, and if there are two, they are necessarily of opposite
sign.
This means that curving tiles do not contribute to $L_t$.
\item
The ribbon crosses itself (Figure \ref{fig:crosssign}).
Consequently there are four crossing of the two boundary strings, two of them
are selfcrossings of one of the strings and do not count, the other two contribute
with the same sign for a total contribution of $\pm 1$ to $L_c$.
The presence of this type of crossings is generally a consequence of the spatial
disposition of the sequence of tiles and in the specific case of face-up planar
configurations (to be considered in Section \ref{sec:piane}) there can be crossings
of this type when we have superposed tiles, or in presence of flap tiles, however the
computation of $L_c$ must be carried out case by case.
\end{enumerate}


\subsection{Contribution of the straight tiles to \texorpdfstring{$L_t$}{Lt}.}
\label{sec:ltstraigth}

The ribbon ``bounces'' exactly once at each straight tile (Definition
\ref{def:localshape}), hence it contributes to $L_t$ with a value $\delta L_t
= \pm \frac{1}{2}$.

After analyzing the various possibilities we conclude for tile $T_i$ as follows:
\begin{itemize}
\item
$\delta L_t = +\frac{1}{2}$ if $T_i$ is a ``vertical'' tile (connected to the
adjacent tiles through its horizontal sides) of type $\sla$, or if it is a horizontal
tile of type $\bsla$;
\item
$\delta L_t = -\frac{1}{2}$ if $T_i$ is a horizontal tile of type $\sla$ or a
vertical tile of type $\bsla$.
\end{itemize}


\subsection{Contribution of the flap tiles to \texorpdfstring{$L_t$}{Lt}.}
\label{sec:ltflaps}

A flap tile can be covered by the ribbon either with four sections (three
``bounces'') of none at all.
In this latter case there is still a ``bounce'' of the ribbon when it
goes from the previous tile to the next (superposed) tile:
the ribbon travels from below the lower tile to above the upper tile or
viceversa.
We need to keep track of this extra bounce.

After analysing the possibilities we conclude for tile $T_i$ as follows:
\begin{itemize}
\item
$\delta L_t = +\frac{1}{2}$ if $T_i$ is a vertical flap of type $\sla$ (connected
to the adjacent tile through a horizontal side), or if it is a horizontal flap
of type $\bsla$;
\item
$\delta L_t = -\frac{1}{2}$ if $T_i$ is a horizontal flap of type $\sla$ or a vertical
flap of type $\bsla$.
\end{itemize}


\subsection{Linking number of constructible configurations.}

\begin{theorem}\label{teo:invariantetopologico}
A constructible spatial configuration of the puzzle necessarily satisfies
$L = 0$.
\end{theorem}

\begin{proof}
The linking number $L$ does not change under legitimate moves of the puzzle, so that
it is sufficient to compute it on the initial configuration of Figure \ref{fig:black2x4}.
There are no superposed tiles nor flaps, so that the ribbon does not cross itself,
hence $L_c = 0$.
The only contribution to $L_t$ comes from the four straight tiles, and using the
analysis of Section \ref{sec:ltstraigth} it turns out that their contribution
cancel one another so that also $L_t = 0$ and we conclude the proof.
\end{proof}


\subsection{Examples of configurations with nonzero linking number.}

Due to Theorem \ref{teo:invariantetopologico} such configurations of the puzzle
cannot be constructed.

\begin{figure}\centering
\includegraphics[width=11cm]{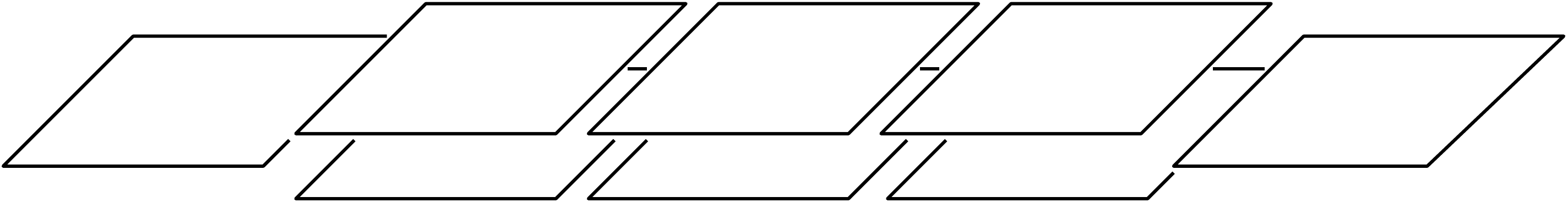}
\caption{\small{This configuration is not constructible because it has linking
number $L \neq 0$.}}
\label{fig:lunga}
\end{figure}

One such configuration is shown in Figure \ref{fig:lunga} and would realize the
maximal possible diameter for a configuration.
The metric invariant of Section \ref{sec:invariantemetrico} is $\Delta = 0$ so
that it is not enough to exclude this configuration,
however we shall show that in this case $L \neq 0$ and conclude that we have a
nonconstructible configuration.
It will be studied in Section \ref{sec:forsecostruibili}.

Another interesting configuration that can be excluded with the topological invariant
and not with the metric one is a ``figure eight'' corresponding to the sequence
\eqref{eq:seqottoobliquo} of Section \ref{sec:forsecostruibili}.
Figure \ref{fig:nc_shapes} (right) show a 3D configuration that cannot be constructed
because $L \neq 0$.

\section{Octominoid configurations}\label{sec:octominoid}
A class of special configurations that can be studied using the two invariants introduced
in Sections \ref{sec:invariantemetrico} and \ref{sec:invariantetopologico} consists of the
so-called \textit{octominoid} configurations.
These correspond to positions of the eight tiles to form a 3D shape with all tiles parallel
to one of the coordinate planes and no pair of superposed tiles.
The term octominoid was introduced by J\"urgen K\"oller in \cite{basteleien:web} and is suggested
by the term \textit{polyominoes} to denote planar shapes made of some fixed number of adjacent
unit squares joined by their sides.

The total number of distinct octominoids is the large number $207265$, most of which can be
immediately excluded as possible configurations of the Rubik's Magic because the eight squares
cannot be cyclically connected by their sides;
then a further reduction is obtained by enforcing the local constraints of Section
\ref{sec:vincolilocali}.

The correct way of counting the set of feasible shapes (configurations of the undecorated
puzzle that satisfy the local constraints) must take into account the cyclical ordering of the tiles
together with the direction of their marked diagonals.
It is thus possible to obtain the same 3D octominoid shape with different puzzle configurations:
they come often in pairs, one configuration obtained from the other by reverting the direction of the diagonals,
but there can be more than two configurations, or just one.

The total number of configurations in the shape of an octominoid that satisfy the local constraints
turns out to be $1291$ realizing a total of $582$ different octominoid shapes.

These numbers are obtained by using a software code that can be downloaded from \cite{mine:web} and
that will be briefly described in Section \ref{sec:code3D}.

Any given octominoid configuration can be described by constructing the so-called \textit{magic code}.
It consists of a sequence of characters (like \texttt{RRmRRmRUmDUm}) with eight capital letters taken
from the set \texttt{RLUD} (standing for \textit{right}, \textit{left}, \textit{up}, \textit{down})
optionally followed by the lower letter \texttt{m} (mountain fold) or \texttt{v} (valley fold).

They encode the relative adjacency information of each of the eight tiles with the next one.
The first tile (say $T_0$) is oriented by drawing an arrow on one of its faces (front) parallel to a side.
The selected orientation for $T_0$ \textbf{must} be such that the nylon strings cross the front face in the direction
south-west to north-east ($\sla$ direction).

The first capital letter indicates which one of the four sides of $T_0$ is connected to $T_1$ (the subsequent tile),
the presence of the lowercase \texttt{m} or \texttt{v} indicates that $T_0$ and $T_1$ form a $90$ degrees angle either
with a mountain fold (letter \texttt{m}) or with a valley fold with respect to the front face of $T_0$.
Otherwise $T_1$ is coplanar with $T_0$.

The orientation of $T_1$ is compatible with the selected orientation of $T_0$ (as defined in
Definition \ref{def:orientation}), i.e. the drawn arrow
on $T_1$ (in case of coplanarity) is exactly the mirror image of the arrow on $T_0$ with respect to the hinging
side (the side of $T_0$ adjacent to $T_1$).
As a simple example the starting $2 \times 4$ configuration of the puzzle can be encoded as \texttt{RRRURRRU} where it
should be noted that the four tiles in the top row have downward arrows.

Since we are interested in configurations of the undecorated puzzle, many distinct magic codes describe the same
configuration based on which tile we select as $T_0$, how we orient it (such that $T_0$ becomes a $\sla$ tile)
and in which direction we traverse the circular chain of eight tiles.
The corresponding magic codes are all considered equivalent.
A canonical magic code is then selected by taking all the disting equivalent magic codes and selecting the first
one with respect to a suitable lexicographic ordering.

The lexicographic comparison is defined such that the four directions are ordered as \texttt{R} $<$ \texttt{U} $<$
\texttt{L} $<$ \texttt{D}, and if the direction is the same, no fold is less than mountain fold which is less than
a valley fold.

The mirror image of a constructible configuration can also be constructed by using the mirror images of the
sequence of construction moves starting from the original $2 \times 4$ configuration, so that we also include all the magic
codes of the mirror images in the same equivalence class.

\begin{figure}
\centering
\includegraphics[width=2.25cm]{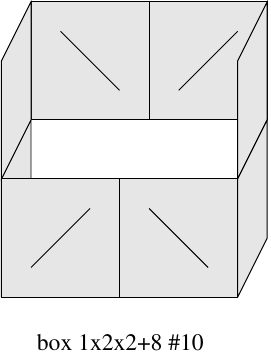}
\includegraphics[width=4.25cm]{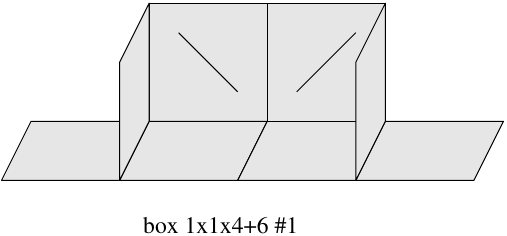}
\includegraphics[width=3.25cm]{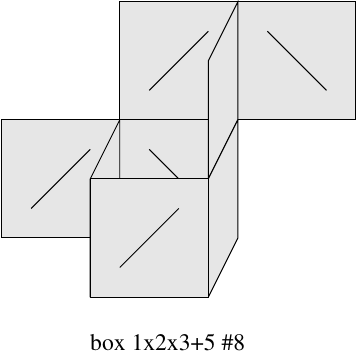}
\caption{\small{Examples of octominoid configurations.  The label indicates the size of the smallest bounding box followed by
$+k$ where $k$ is the number of tiles that lie completely in the boundary of the bounding box.
Finally $\#n$ is the sequential number in the table of the octominoid shapes in \cite[3D octominoid
shapes]{mine:web}.}}\label{fig:shapes}
\end{figure}

\begin{figure}
\centering
\includegraphics[height=9.00cm]{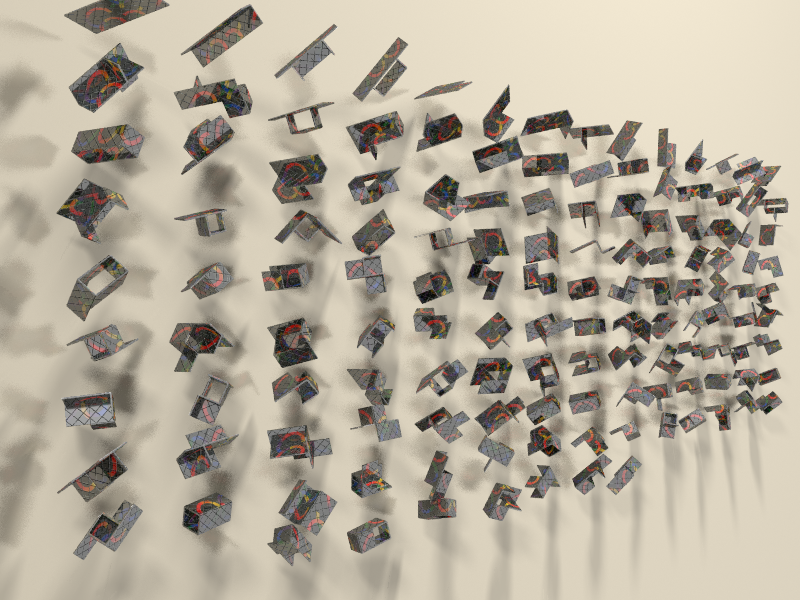}
\includegraphics[height=9.00cm]{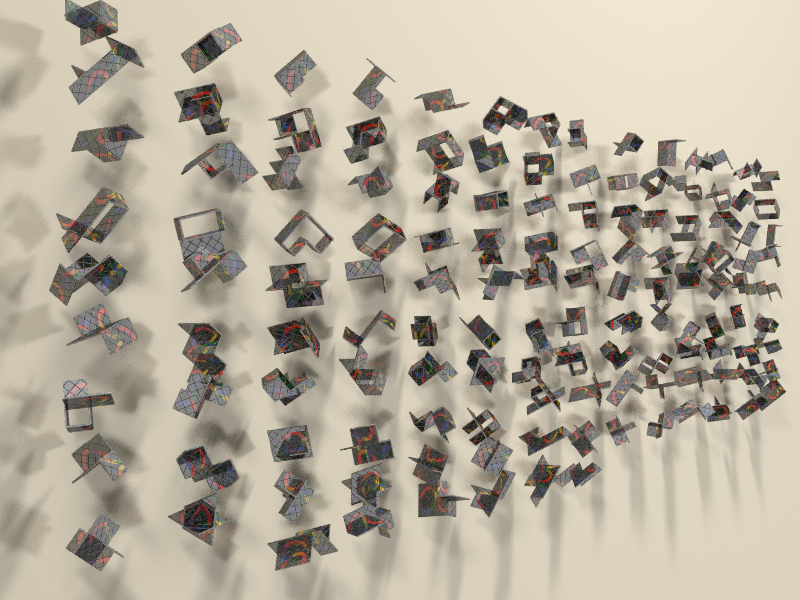}
\caption{\small{All $265$ octominoid shapes with vanishing invariants.
Syntetic images obtained with PovRay.}}\label{fig:tutti}
\end{figure}

Figure \ref{fig:shapes} shows three examples of octominoid configurations, the shape on the left is encoded
by the magic code \texttt{RRmULmLLmURm}, obtained by selecting as $T_0$ the lower-left front tile, oriented with
an upward arrow in the visible face and traversing the configuration counterclockwise.
Then tile $T_1$ is hinged at the right side of $T_0$ (hence the first \texttt{R} character in the magic code) and is
coplanar with $T_0$ (no \texttt{m} or \texttt{v} character following the first \texttt{R}.
Tile $T_2$ is also hinged at the right side of $T_1$ and tilted $90$ degrees with a mountain fold (\texttt{Rm}).
Tile $T_3$ is hinged at the upper side of $T_2$, note that the orientation of $T_3$ is consequently oriented with
a downward arrow.

We have a number of equivalent magic codes by changing the choice and orientation of the starting tile, however the
one coming first in the lexicographic ordering (the canonical magic code) is the string described above.

Similarly, the canonical magic code of the second image of Figure \ref{fig:shapes} is given by \texttt{RRRLmDvUDvLm}, obtained
by selecting as $T_0$ the leftmost tile.
The third configuration of the figure is finally encoded by \texttt{RUUmLDvUvDmD} by suitably chosing and orienting
the tile $T_0$.

In the special class of octominoids, we were able to automate the computation of both the metric and
the topological invariants of the configuration having a given magic code, thus allowing to quickly rule out all
configurations with nonzero invariants.
Enforcing zero metric invariant reduces the configurations from $1291$ down to $737$, then enforcing also the vanishing of
the linking number further reduces the number of configurations to $460$ ($265$ distinct octominoid shapes).

These are all listed in \cite[3D octominoids]{mine:web}, ordered according to the dimension of the smallest bounding
box.
Surprisingly, all of them could be actually constructed with the real puzzle, and building istructions are available
in \cite[3D octominoids]{mine:web}.
A few of such shapes are not present in the comprehensive table of symmetric 3D shapes in
\cite{basteleien:web} and might be possibly obtained by us for the first time.

Figure \ref{fig:tutti} shows all $265$ such shapes in one shot.
The images where obtained synthetically using the ray-tracing software
\texttt{PovRay}, in conjunction with the output of our software code to
obtain the list of admissible magic codes.


\section{Planar face-up configurations}\label{sec:piane}
We shall apply the results of the previous sections to a particular choice
of spatial configurations, we shall restrict to planar configurations
(all tiles parallel to the horizontal plane) with non-overlapping consecutive
tiles.
Superposed nonconsecutive tiles are allowed.

They can be obtained starting from strings of cardinal directions in the following way.

The infinite string
$s : \Z \to \{E,N,W,S\}$
is a typographical sequence with index taking values in the integers $\Z$ where the
four symbols stand for the four cardinal directions East, North, West, South.
On $s$ we require

\begin{enumerate}
\item
Periodicity of period $8$: $s_{n+8} = s_n$ for any $n \in \Z$;
\item
Zero mean value: in any subsequence of $8$ consecutive characters (for example in
$\{s_0, \dots, s_7\}$) there is an equal number of characters $N$ as of characters
$S$ and of characters $E$ as of characters $W$.
\end{enumerate}

An \textbf{admissible sequence} is one that satisfies the two above requirements.

Periodicity allows us to describe an admissible sequence by listing $8$ consecutive
symbols, for definiteness and simplicity we shall then describe an admissible sequence
just by listing the symbols $s_1$ to $s_8$.


The character $s_i$ of the string indicates the relative position between the
two consecutive tiles $T_{i-1}$ and $T_i$, that are horizontal and face-up.

The first tile $T_0$ can be of type $\sla$ or type $\bsla$, all the others $T_i$ are
of the same type as $T_0$ if $i$ is even, of the opposite type if $i$ is odd.
The local constraints allows to recover a spatial configuration of the puzzle from
an admissible sequence with two caveats:

\begin{enumerate}
\item
For at least one of the tiles, say $T_0$, it is necessary to specify if it is of type
$\sla$ or $\bsla$.
We can add this information by inserting the symbol $\sla$ or $\bsla$ between two
consecutive symbols, usually before $s_1$;
\item
In case of superposed tiles (same physical position) it is necessary to clarify
their relative position (which is above which).
We can add this information by inserting a positive natural number between
two consecutive symbols that indicates the ``height'' of the corresponding tile.
In the real puzzle the tiles are not of zero width, so that their height in space
cannot be the same.
In case of necessity we shall insert such numbers as an index of the symbol at the
left.
\end{enumerate}

\begin{remark}[Configurations that can be assembled]\label{rem:assemblabili}
Given an admissible sequence it is possible to compute the number of superposed tiles
at any given position.
An \textbf{assemblage} of a sequence entails a choice of the height of each of the
superposed tiles (if there is more than one).
We do this by adding an index between two consecutive symbols.
However for this assemblage to correspond to a possible puzzle configuration we need
to require a condition.
We hence say that an assemblage is \textbf{admissible} if whenever tile $T_i$ is
superposed to tile $T_j$, $i \neq j$, and also $T_{i \pm 1}$ is superposed to
$T_{j+1}$, then the relative position of the tiles in the two pairs cannot be exchanged.
This means that if $T_i$ is at a higher height than $T_j$, then $T_{i \pm 1}$ cannot be
at a lower height than $T_{j+1}$.
It is possible that a given admissible sequence does not allow for any admissible
assemblage or that it can allow for more than one admissible assemblage.
\end{remark}

Observe that the mirror image of an oriented spatial configuration of the undecorated
puzzle entails a change of type, $\sla$ tiles become of type $\bsla$ and viceversa.
If the mirror is horizontal the reflected image is a different assemblage of the same
admissible sequence with all tiles of changed type and an inverted relative position
of the superposed tiles.

On the set of admissible sequences we introduce an equivalence relation defined by
$s \equiv t$ if one of the following properties (or a combination of them) holds:

\begin{enumerate}
\item (cyclicity) The two sequences coincide up to a translation of the index:
$s_n = t_{n+k}$ for all $n$ and some $k \in \Z$;
\item (order reversal) $s_n = t_{k-n}$ for all $n$ and some $k \in \Z$;
\item (rotation) $s$ can be obtained from $t$ after substituting $E \to N$, $N \to W$,
$W \to S$, $S \to E$;
\item (reflection) $s$ can be obtained from $t$ after substituting $E \to W$ e $W \to E$.
\end{enumerate}

Let us denote by $\SSS$ the set of equivalence classes.

We developed a software code capable of finding a canonical representative of
each of these equivalence classes, they are $71$ (cardinality of $\SSS$).
In Table \ref{tab:sequenze} we summarize important properties of these canonical
sequences, subdivided with respect to the number of flap tiles.
It is worth noting that some of the $71$ sequences admit more than one nonequivalent
admissible assemblages in space due to the arbitrariness in choosing the type of tile
$T_0$ and the ordering of the superposed tiles.
A few of the $71$ admissible sequences do not admin any admissible assemblage, one
of these is the only sequence with $8$ flaps: $EWEWEWEW$.
Since the constructability of a spatial configuration is invariant under specular
reflection (which entails a change of type of all tiles) we can fix the type of tile
$T_0$, possibly reverting the order of the superposed tiles.

The canonical representative of an equivalence class in $\SSS$ is selected by
introducing a lexicographic ordering in the finite sequence $s_1, \dots, s_8$ where
the ordering of the four cardinal directions is fixed as $E < N < W < S$.
Then the canonical representative is the smallest element of the class with respect to
this ordering.

The source of the software code can be downloaded from the web page \cite{mine:web}.


\begin{table}[ht]
\caption{\small{Sequences in $\SSS$.}}
\centering
\begin{tabular}{c c c c c c c}
\hline\hline
number & \small{number of} & \small{number of}
                                & $\Delta = 0$ & $\Delta = L = 0$ & not \\
of flaps & sequences & assembl. & assembl. & assembl. & classified \\
\hline
none & 7 & 6 & 4 & 2 &  - \\
1       & 7 & 14 & 7 & 3 & - \\
2       & 22 & 44 & 20 & 13 & 3 \\
3       & 10 & 50 & 15 & 5 & 2 \\
4       & 18 & 38 & 11 & 2 & - \\
5       & 2 & 12 & 1 & - & - \\
6       & 4 & 4 & 1 & - & - \\
8       & 1 & 0 & - & - & - \\
\hline
total  & 71 & 168 & 59 & 25 & 5 \\
\end{tabular}
\label{tab:sequenze}
\end{table}


In Table \ref{tab:2e4flap} the sequences with two and four flaps are subdivided
based on the distribution of the flaps in the sequence.

\begin{theorem}\label{teo:ltiszero}
All planar face-up configurations have zero ``twist'' contribution to the
topological invariant: $L_t = 0$.
Consequently we have $L = L_c$ and to compute the linking number it is sufficient to
compute the contributions coming from the crossing of the ribbon with itself.
Any planar face-up configuration with an odd number of selfintersectons of the
ribbon with itself has $L \neq 0$.
\end{theorem}

\begin{proof}\footnote{%
This proof is due to Giovanni Paolini, Scuola Normale
Superiore of Pisa.}
We denote with $k_1$, ..., $k_s$ the number of symbols in contiguous subsequence of
$E$, $W$ (horizontal portions) or of $N$, $S$ (vertical portions).
Each portion of $k_i$ symbols contains $k_i - 1$ straight tiles or flaps, all
``horizontal'' or ``vertical'', hence each tile contributes to $L_t$ with alternating
sign due to the fact that the tiles are alternatively of type $\sla$ and $\bsla$.
If $k_i-1$ is even, then the contribution of this portion is zero, while if it is
odd it will be equal to the contribution of the first straight or flap tile of the
portion.
It is not restrictive to assume that the first portion of $k_1$ symbols is horizontal
and the last (of $k_s$ symbols) is vertical.
In this way if $i$ is odd, then $k_i$ is the number of symbols in a horizontal
portion whereas if $i$ is even, then $k_i$ is the number of symbols in a vertical
portion.
Up to a change of sign of $L_t$ we can also assume that the first tile is of type
$\sla$.
Finally we observe that $k_i > 0$ for all $i$.
Twice the contribution to $L_t$ of the $i$-th portion is given by
\begin{equation}\label{eq:ltparziale}
(-1)^{i-1} (-1)^{k_1 + k_2 + ... + k_{i-1}} (1 + (-1)^{k_i})
\end{equation}
where the last factor in parentheses is zero if $k_i$ is odd and is $2$ if $k_i$
is even;
the sign changes on vertical portions with respect to horizontal portions
(factor $(-1)^{i-1}$) and changes when the type ($\sla$ or $\bsla$) of the first
straight or flap til of the portion changes (factor $(-1)^{k_1 + k_2 + ... + k_{i-1}}$).
Summing up \ref{eq:ltparziale} on $i$ and expanding we have
\begin{eqnarray*}
2 L_t &= & \displaystyle{\sum_{i=1}^s (-1)^{i-1} (-1)^{k_1 + k_2 + ... + k_{i-1}}
+
\sum_{i=1}^s (-1)^{i-1} (-1)^{k_1 + k_2 + ... + k_{i-1}} (-1)^{k_i}}
\\
&= & \displaystyle{- \sum_{i=1}^s (-1)^i (-1)^{k_1 + k_2 + ... + k_{i-1}}
+
\sum_{i=1}^s (-1)^{i+1} (-1)^{k_1 + k_2 + ... + k_i}}
\\
&= & \displaystyle{- \sum_{i=1}^s (-1)^i (-1)^{k_1 + k_2 + ... + k_{i-1}}
+
\sum_{i=2}^{s+1} (-1)^i (-1)^{k_1 + k_2 + ... + k_{i-1}}}
\\
&= & 1
+
 (-1)^{s+1} (-1)^{k_1 + k_2 + ... + k_s}
= 0
\end{eqnarray*}
because $s$ is even and $k_1 + \dots + k_s = 8$, even.
\end{proof}

\section{Configurations with vanishing invariants}\label{sec:forsecostruibili}
We shall identify admissible assemblages whenever they correspond to equivalent
puzzle configurations, where we also allow for specular images.
In particular this allows us to assume the first tile to be of type $\sla$.

Assemblages corresponding to non-equivalent sequences cannot be equivalent,
on the contrary there can exist equivalent assemblages of the same sequence and
this typically happens for symmetric sequences.

The two invariants can change sign on equivalent sequences or equivalent assemblages,
this is not a problem since we are interested in whether the invariants are zero or
nonzero.
In any case the computations are always performed on the canonical representative.

The contribution $\Delta_c$ of $\Delta = \Delta_c + \Delta_f$ (coming from the
curving tiles) can be computed on the sequence (it does not depend on the assemblage).
On the contrary the contribution $\Delta_f$ coming from the flap tiles depends on
the actual assemblage.

With the aid of the software code we can partially analyze each canonical
admissible sequence and each of the possible admissible assemblages of a sequence.
In particular the software is able to compute the metric invariant of an
assemblage, so that we are left with the analysis of the topological invariant, and
we shall perform such analysis only on assemblages having $\Delta = 0$, since
our aim is to identify as best as we can the set of constructible configurations.


\subsection{Sequences with no \textit{flaps}}\label{sec:noflaps}


There are seven such sequences, three of them do not have any superposed tiles,
so that they cover a region of the plane corresponding to $8$ tiles (configurations
of area $8$).
For these three sequences we only have one possible assemblage (having fixed the
type $\sla$ of tile $T_0$).

\noindent%
The sequence
\begin{equation}\label{eq:seq2x4}
EEENWWWS
\end{equation}
corresponds to the initial configuration $2\times 4$ of the puzzle.
The sequence
\begin{equation}\label{eq:seq3x3}
EENNWWSS
\end{equation}
corresponds to the ``window shape'', a $3\times 3$ square without the central tile.
The sequence
\begin{equation}\label{eq:seqelle}
EENNWSWS
\end{equation}
corresponds to the target configuration of the puzzle
(Figure \ref{fig:blacksolved}).
Two sequences cover $7$ squares of the plane (area $7$),
the sequence
\begin{equation}\label{eq:seqottoobliquo}
EENWSSWN
\end{equation}
and the sequence
\begin{equation}\label{eq:seqbirillo}
ENENWSWS .
\end{equation}
A sequence without flaps and area $6$ (two pairs of superposed tiles) is
\begin{equation}\label{eq:seqottoverticale}
ENESWNWS .
\end{equation}
The last possible sequence (with area $4$) would be
\begin{equation}\label{eq:seqspirale}
ENWSENWS ,
\end{equation}
this however cannot be assembled in space since it consists of a closed circuit of
$4$ tiles traveled twice
(see Remark \ref{rem:assemblabili}).

The metric invariant is nonzero (hence the corresponding assemblage is not constructible) for the two sequences \eqref{eq:seq3x3} and \eqref{eq:seqbirillo},
the topological invariant $L$ further reduces the number of possibly constructible
configuration by excluding also the two sequences
\eqref{eq:seqottoobliquo} e \eqref{eq:seqottoverticale}.

The remaining two configurations, corresponding to sequences
\eqref{eq:seq2x4} ed \eqref{eq:seqelle}, are actually constructible
(Figures \ref{fig:black2x4} and \ref{fig:blacksolved}).


\subsection{Sequences with one \textit{flap}}

We find seven (nonequivalent) sequences with exactly one flap.
Three of these have area $7$:
\begin{equation}\label{eq:seqf1a7a}
EENNWSSW   
\end{equation}
\begin{equation}\label{eq:seqf1a7b}
EENWNSWS
\end{equation}
\begin{equation}\label{eq:seqf1a7c}
EEENWWSW   
\end{equation}
and four have area $6$:
\begin{equation}\label{eq:seqf1a6a}
EENWSWSN 
\end{equation}
\begin{equation}\label{eq:seqf1a6b}
EENWSWNS
\end{equation}
\begin{equation}\label{eq:seqf1a6c}
EENWSNWS
\end{equation}
\begin{equation}\label{eq:seqf1a6d}
EEENWSWW .  
\end{equation}

In all cases it turns out that there are two nonequivalent assemblages of each
of these sequences according to the flap tile being ascending or descending, and
they have necessarily a different value of $\Delta$, so that at most one (it turns
out exactly one) has $\Delta = 0$.
We shall restrict the analysis of the topological invariant to those having $\Delta
= 0$.

The two sequences \eqref{eq:seqf1a7a} and \eqref{eq:seqf1a7c} have $\Delta = 0$
if the (horizontal) flap tile is descending (Figure \ref{fig:nastroflap} right).
The linking number reduces to $L = L_c$ (Theorem \ref{teo:ltiszero}).
Since in both cases we have exactly one crossing of the ribbon with itself
we conclude that $L \neq 0$ and the sequences are \textbf{not constructible}.

To have $\Delta = 0$ the vertical flap of the sequence \eqref{eq:seqf1a7b} must be
descending.
Then there is one crossing of the ribbon with itself, so that 
$L \neq 0$ and the configuration is \textbf{not constructible}.

Sequences \eqref{eq:seqf1a6a} and \eqref{eq:seqf1a6b} have $\Delta = 0$ provided their
flap is ascending.
We have now two crossings of the ribbon with itself and they turn out to have
opposite sign in their contribution to $L_c$, so that $L = 0$ and the two sequences
``might'' be constructible.

Sequences \eqref{eq:seqf1a6c} and \eqref{eq:seqf1a6d} have $\Delta = 0$ provided their
flap is ascending.
Sequence \eqref{eq:seqf1a6c} is then \textbf{not constructible} because there is
exactly one selfcrossing of the ribbon so that $L = L_c \neq 0$.
On the contrary, sequence \eqref{eq:seqf1a6d} exhibits two selfcrossings with opposite
sign and $L = L_c = 0$.

In conclusion of the $7$ different sequences with one flap, four are necessarily
non constructible because the topological invariant is non-zero, the remaining three
sequences: \eqref{eq:seqf1a6a},
\eqref{eq:seqf1a6b},
\eqref{eq:seqf1a6d}
are actually constructible (refer to \cite[planar face-up]{mine:web} for building instructions:
first, second and third image of section ``One flap"), see also Figure \ref{fig:constructible1},
first three images.


\subsection{The two sequences with two adjacent \textit{flaps}}

%
%
Adjacency of the two flaps entails that both are ascending or both descending
(Remark \ref{rem:assemblabili}) and also they are both horizontal or both vertical
since they are hinged to each other so that they contribute to the metric invariant
$\Delta_f = \pm 4$ whereas $\Delta_c = 0$.
Hence the metric invariant is nonzero and the two sequences are non-constructible.


\subsection{The five sequences with two \textit{flaps} separated by one tile}

\subsubsection{Sequence $EENWSEWW$}
$\Delta = 0$ implies that the two (horizontal) flaps are one ascending and one
descending.
There are two non-equivalent admissible assemblages satisfying $\Delta = 0$,
computation of the topological invariant gives $L = L_c = \pm 4$ for one of the two
assemblages whereas the other has $L = L_c = 0$ and is actually constructible
(Figure \ref{fig:constructible2}, second image):
\begin{equation}\label{eq:seqf2a5c}
\sla E_3 E_2 N W S_2 E_1 W_1 W .
\end{equation}
For building instructions, see \cite[planar face-up]{mine:web}, second image of section
``Two flaps, area $5$''.
%
%

\subsubsection{Sequence $ENEWSNWS$}
$\Delta = 0$ implies that both flaps (one is horizontal and one vertical) are
ascending or both descending.
The two corresponding distinct admissible assemblages have both $L = L_c  = 0$.
The two assemblages are:
\begin{equation}\label{eq:seqf2a5ab}
\sla E_2 N_3 E W_2 S_2 N_3 W S
\quad \text{and} \quad
\sla E_1 N_1 E W_2 S_2 N_3 W S .
\end{equation}
They are both constructible (Figure \ref{fig:constructible2}, third and first image
respectively).
For building instructions: \cite[planar face-up]{mine:web}, third and first image
of section ``Two flaps, area $5$''.
%
%
%
%

\subsubsection{Sequence $ENEWNSWS$}
We can fix the first tile $T_0$ to be of type $\sla$, then $\Delta_c = 4$ and
$\Delta = 0$ implies that the first flap (horizontal) is ascending and the second
(vertical) is descending.
Computation of the topological invariant gives $L = L_c = 0$ and we have another
unclassified sequence:
\begin{equation*}
\sla E N_1 E W_3 N S_2 W S .
\end{equation*}
The two lowest superposed tiles can be exchanged, however the resulting assemblage
is equivalent due to the reflection symmetry of the sequence of symbols.

%

\subsubsection{Sequence $ENWESNWS$}
Imposing $\Delta = 0$ the two flaps (one is horizontal and one is vertical) must
be both ascending or both descending.
In both cases we compute $L = L_c = 0$.
Actually the two assemblages are equivalent by taking advantage of the symmetry
of the sequence, one of these is
\begin{equation}\label{eq:seqf2a4a}
\sla E_1 N_1 W_1 E_2 S_2 N_3 W_2 S
\end{equation}
and is constructible (Figure \ref{fig:constructible2}, fourth image);
building instructions in \cite[planar face-up]{mine:web}, image in section ``Two flaps, area $4$''.

%

\subsubsection{Sequence $EENWSSNW$}
Imposing $\Delta = 0$ the two flaps (one horizontal and one vertical) must
be both ascending or both descending.
In both cases we compute $L = L_c = 0$.
The two assemblages are equivalent as in the previous case, one of these
(Figure \ref{fig:constructible1}, $8$-th image) is
\begin{equation}\label{eq:seqf2a6d}
\sla E_1 E N W S_3 S N_2 W 
\end{equation}
and building instructions can be found in \cite[planar face-up]{mine:web}, fifth image of
section ``Two flaps, area $6$.

%
%

\subsection{The three sequences with two \textit{flaps} separated by two tiles}

\begin{figure}\begin{center}
\includegraphics[width=12cm]{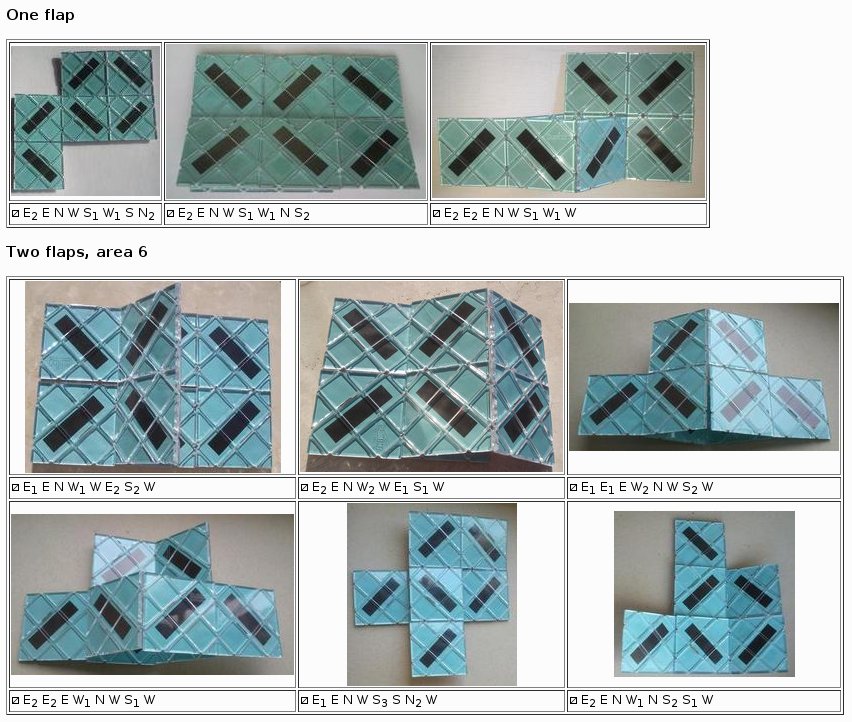}
\caption{\small{The three constructible configurations with one flap and those with two flaps and
area $6$.
For building instructions we refer to \cite{mine:web}, click on ``planar face-up''.}}\label{fig:constructible1}
\end{center}\end{figure}

%
%
All three admissible sequences with two flaps at distance $3$ (separated by
two tiles) have $\Delta_c = 0$.
Two of these sequences have both horizontal or both vertical flaps, so that
$\Delta = 0$ entails that one flap is ascending and one is descending.
The third sequence has an horizontal flap and a vertical flap so that $\Delta = 0$
entails that both flaps are ascending or both descending.
In all cases we have two selfcrossings of the ribbon with opposite sign, hence
$L = L_c = 0$ and might be constructible.
Each of the three sequences admit two distinct assemblages both with $\Delta
= L = 0$:
\begin{eqnarray}
\sla E_2 E_2 E W_1 N W S_1 W
\quad \text{,} \quad
\sla E_1 E_1 E W_2 N W S_2 W \label{eq:seqf2a6bc}
\\
\sla E_2 E N W_1 N S_2 S_1 W
\quad \text{,} \quad
\sla E_1 E N W_2 N S_1 S_2 W \label{eq:seqf2a6e}
\\
\sla E_2 E N W_2 W E_1 S_1 W
\quad \text{,} \quad
\sla E_1 E N W_1 W E_2 S_2 W \label{eq:seqf2a6a}
\end{eqnarray}
The first two and the last two are actually constructible
(Figure \ref{fig:constructible1}, $4$-th, $5$-th, $6$-th and $7$-th image respectively);
building instructions in \cite[planar face-up]{mine:web}, images $1$ to $4$ of section
``Two flaps, area $6$''.
The first of the second row is also constructible (Figure \ref{fig:constructible1},
$9$-th image)
although with a considerable amount of strain on the nylon wires.
Building instructions in \cite[planar face-up]{mine:web}, $6$-th image of
section ``Two flaps, area $6$''.

\begin{figure}\begin{center}
\includegraphics[width=12cm]{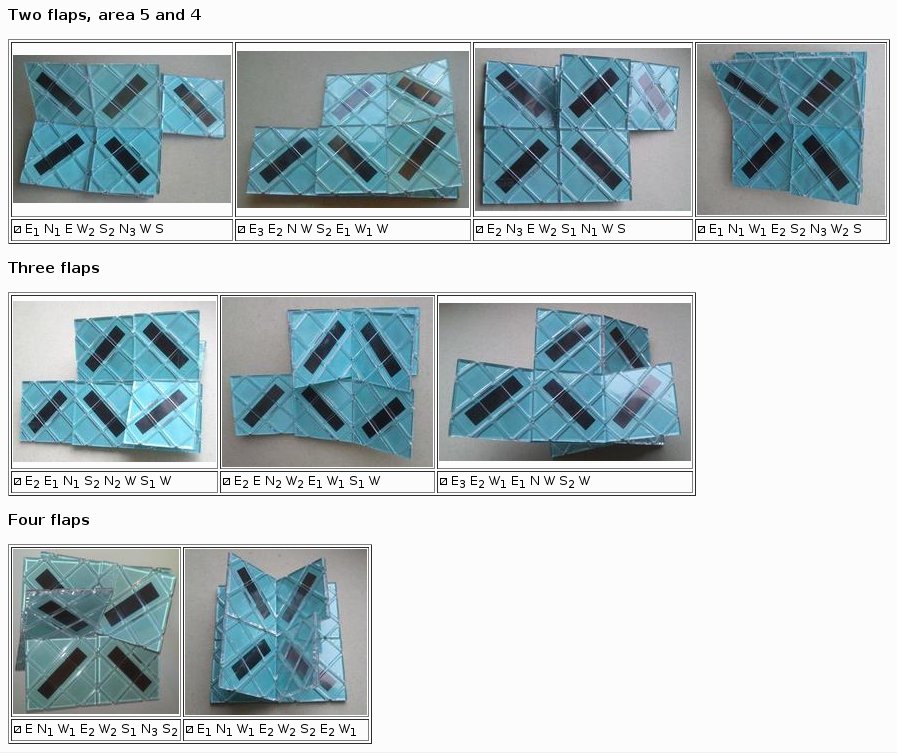}
\caption{\small{Constructible configurations with two flaps, area $5$ and $4$ and configurations with
three and four flaps.
For building instructions we refer to \cite{mine:web}, click on ``planar face-up''.}}\label{fig:constructible2}
\end{center}\end{figure}

\subsection{The twelve sequences with two flaps in antipodal position}

\begin{table}[ht]
\caption{\small{Sequences with two flaps (left) and four flaps (right) subdivided
based on the relative position of the flaps.}}
\centering
\begin{tabular}{c c}
\hline\hline
sequences & distribution \\
with $2$ flaps & of flaps \\
\hline
2 & ffxxxxxx \\
5 & fxfxxxxx \\
3 & fxxfxxxx \\
12 & fxxxfxxx \\
\hline
\end{tabular}
~~~
\begin{tabular}{c c}
\hline\hline
sequences & dist. of flaps \\
\hline
1 & ffffxxxx \\
5 & fffxxfxx \\
1 & ffxffxxx \\
4 & ffxfxxfx \\
1 & ffxxffxx \\
6 & fxfxfxfx \\
\hline
\end{tabular}
\label{tab:2e4flap}
\end{table}

%

Of the $12$ sequences with two flaps in opposite (antipodal) position we
first analyze those (they are $10$) in which the tiles follow the same path
from one flap to the other and back.
One of these is shown in Figure \ref{fig:lunga}.
All have $\Delta_c = 0$ so that the contribution of the two flaps must have
opposite sign in order to have $\Delta = 0$.
If one flap is horizontal and the other vertical, then they must be both
ascending or both descending and we have no possible admissible assemblage
(Remark \ref{rem:assemblabili}).
We are then left with those sequences having both horizontal or both vertical
flaps, one ascending and one descending.
In this situation we find that the ribbon has $3$ selfcrossings, so that
necessarily $L = L_c \neq 0$ and these sequences are also not constructible.

We remain with the two sequences $EENEWWSW$ and $EENNSWSW$ that both have a
contribution $\Delta_c = -4$ (fixing $T_0$ of type $\sla$),
so that the two flaps must contribute with a positive sign to the metric invariant.
The first sequence has both horizontal flaps, and they must be both descending,
this is now possible thanks to the different path between the two flaps.
The second sequence has one horizontal and one vertical flap, so that the first
must be ascending and the second descending.
There are exactly two selfcrossings of the ribbon in both cases, however they
have the same sign in the first case implying $L = L_c \neq 0$, hence
non constructible.
They have opposite sign in the second case and we have both zero invariants.
In conclusion the only one of the $12$ sequences that might be constructible
is
\begin{equation}
\sla E_1 E N_1 N S_2 W S_2 W .
\end{equation}


\subsection{Sequences with three \textit{flaps}}

%
%
%

Of the $10$ admissible sequences with three flaps there is only one
with all adjacent flaps, having $\Delta_c = \pm 2$.
The three flaps being consecutive are all horizontal or all vertical and
all ascending or all descending, with a total of $\Delta_f = \pm 6$ and
the metric invariant cannot be zero.

%
%
%
%
Four sequences have two adjacent flaps, and in all cases $\Delta_c = \pm 2$.
Imposing $\Delta = 0$ allows to identify a unique assemblage for each sequence
(with one exception).
In all cases a direct check allows to compute $L = L_c = 0$.
These sequences are:
\begin{eqnarray}
\sla E_3 E_2 W_{1,2} E_1 N W S_{2,1} W
\quad \text{,} \quad
\sla E_3 E N W_1 S_1 N_2 S_2 W \label{eq:seqf3a5c}
\\
\sla E_2 E N_2 W_2 E_1 W_1 S_1 W
\quad \text{,} \quad
\sla E_2 E_1 N_1 S_2 N_2 W S_1 W . \label{eq:seqf3a5a}
\end{eqnarray}
The last two are actually constructible (Figure \ref{fig:constructible2}, $5$-th and $6$-th images),
building instructions in \cite[planar face-up]{mine:web}, first and second images of section
``Three flaps''.
One of the two assemblages of the left sequence in \eqref{eq:seqf3a5c} can be
actually constructed (Figure \ref{fig:constructible2}, $7$-th image),
building instructions in \cite[planar face-up]{mine:web}, third image of section
``Three flaps''.
If the puzzle had sufficiently deformable nylon threads and tiles we could
conceivably deform the first assemblage into the second.
We do not know at present if the right sequence in \eqref{eq:seqf3a5c} is
constructible (unclassified).
%
%
%
%

The five remaining sequences all have $\Delta_c = \pm 2$.
Imposing $\Delta = 0$ leaves us with $10$ different assemblages:
the sequence with all three horizontal flaps has three different assemblages
with $\Delta = 0$,
three of the remaining four sequences (with two flaps in one direction and the
third in the other direction) have two assemblages each,
the remaining sequence has only one assemblage with $\Delta = 0$.
In all cases a direct check quantifies in $3$ or $5$ (in any case an odd value)
the number of selfcrossings of the ribbon, so that $L = L_c \neq 0$.
None of these sequences is then constructible.


\subsection{Sequences with four \textit{flaps}}

There are $18$ such sequences.
Six of these have a series of at least three consecutive flaps and a contribution
$\Delta_c = 0$.
They are not constructible because the consecutive flaps all contribute with the
same sign to $\Delta_f$.



The sequences $ENWEWSNS$ and $ENWEWSEW$ have $\Delta_c = 0$ and two pairs of
adjacent flaps oriented in different directions in the first case and in the same
direction in the second case.
To have $\Delta = 0$ they must contribute with opposite sign and hence must be
all ascending or all descending in the first case whereas in the second case
one pair of flaps must be ascending and one descending.
Thanks to the symmetry of the sequences the two possible assemblages of each are
actually equivalent.
The linking number turns out to be $L = 0$ and we have two possibly constructible
configurations:
\begin{equation}
\sla E N_1 W_1 E_2 W_2 S_1 N_3 S_2
\quad \text{and} \quad
\sla E_1 N_1 W_1 E_2 W_2 S_2 E_2 W_1 .
\label{eq:seqf4a4ab}
\end{equation}
Both turn out to be constructible (Figure \ref{fig:constructible2}, $8$-th and $9$-th images),
building instructions in \cite[planar face-up]{mine:web}, images of section ``Four flaps''.

There are four sequences with a single pair of adjacent flaps, the other two
being isolated, all with $\Delta_c = 0$.
The two isolated flaps must contribute to the metric invariant with the same sign,
opposite to the contribution that comes from the two adjacent flaps.
In three of the four cases the two isolated flaps have the same direction and
hence both must be ascending or both descending.
It turns out that there is no admissible assemblage with such characteristics.
The pair of adjacent flaps of the remaining sequence ($EENSWEWW$) are horizontal.
If they are ascending the remaining horizontal flap must be descending whereas
the vertical flap must be ascending (to have $\Delta = 0$).
This situation (or the one with a descending pair of adjacent flaps) is assemblable
and we can compute the linking number, which turns out to be
$L = L_c = \pm 2$.
Even this configuration is not constructible.

The remaining six sequences (flaps alternating with non-flap tiles) all have
the non-flap tiles superposed to each other.
%
%
An involved reasoning, or the use of the software code, allows to show that for
two of this six sequences, having area $3$, namely
$EEWWEEWW$ and $ENSWENSW$, there is no possible admissible assemblage with
$\Delta = 0$.

\subsubsection{Sequence \texorpdfstring{$ENSEWNSW$}{ENSEWNSW}}
This sequence has area $4$, with two of the four flaps superposed to each
other.
Using the software code we find two different assemblages having $\Delta = 0$.
%
%
Computation of the linking number leads in both cases to
$L = L_c = \pm 2$, hence this sequence is not constructible.

\subsubsection{Sequence \texorpdfstring{$EEWNSEWW$}{EEWNSEWW}}
This sequence has also area $4$, with two of the four flaps superposed to each
other.
Using the software code we find only one assemblages having $\Delta = 0$.
%
%

%
Computation of the linking number leads in to
$L = L_c = \pm 2$, hence this sequence is also not constructible.

\subsubsection{Sequence \texorpdfstring{$EEWNSSNW$}{EEWNSSNW}}
This sequence has area $5$ with no superposed flaps and contribution
$\Delta_c = 0$, so that to have $\Delta = 0$ two flaps contribute positively
and two contribute negatively to the metric invariant.

The software code gives three different assemblages with $\Delta = 0$.

An accurate analysis of the selfcrossings of the ribbon due to the flaps
shows that flaps that contribute positively to the metric invariant also
contribute with an odd number of selfcrossings of the ribbon, besides there
is one selfcrossing due to the crossing straight tiles.

In conclusion we have an odd number of selfcrossings, hence the sequence
is not constructible.

\subsubsection{Sequence \texorpdfstring{$ENSEWSNW$}{ENSEWSNW}}
This sequence has also area $5$ with no superposed flaps, but now the contribution
of the curving tiles to the metric invariant is
$\Delta_c = 4$, so that to have $\Delta = 0$ exactly one of the flaps has positive
contribution to the metric invariant.
This flap will also contribute with an odd number of selfcrossings of the ribbon.
In this case there are no other selfcrossings of the ribbon because there are no
straight tiles, so that we again conclude that the number of selfcrossings of the
ribbon is odd and that $L = L_c \neq 0$.
This sequence is also not constructible.


\subsection{Sequences with five \textit{flaps}}

Both sequences have $\Delta_c = -2$.


The software code quickly shows that the sequence $EEWNSNSW$ does not have any
admissible assemblage with $\Delta = 0$.

The other sequence is $EEWEWNSW$ and to have $\Delta = 0$ the three consecutive
horizontal flaps must be ascending, and also the vertical flap must be ascending.
There is one admissible assemblage satisfying these requirements, but the
resulting number of selfcrossings of the ribbon is odd, and so also this
configuration is not constructible.

%
%


\subsection{Sequences with six \textit{flaps}}

None of the four sequences with six flaps is constructible.
Indeed it turns out that all have $\Delta_c = 0$, so that in order to have
$\Delta = 0$ they must have three flaps with positive contribution and three
with negative contribution to $\Delta_f$.

Two of the four sequences have four or more flaps that are consecutive and
hence all contribute with the same sign to $\Delta_f$, so that $\Delta \neq 0$.

The sequence $EEWEWWEW$ must have three ascending consecutive flaps and three
consecutive descending flaps (all flaps are horizontal).
Analyzing the ribbon configuration shows that there are an odd number of ribbon
selfcrossings, hence $L = L_c \neq 0$.

Finally, all the six flaps of sequence $ENSNSWEW$ must be ascending in order
to have $\Delta = 0$, there is no admissible assemblage with this property.



\subsection{Sequences with seven \textit{flaps}}

There is none.


\subsection{Sequences with eight \textit{flaps}}

The only one is $EWEWEWEW$, but there is no admissible assemblage of this
sequence.


\section{Unclassified configurations}\label{sec:nonclassificate}
\subsubsection{Sequences with two \textit{flaps}}

There are $3$ unclassified assemblages, corresponding to $3$ sequences.
The assemblage
\begin{equation*}
\sla E N_1 E W_3 N S_2 W S
\end{equation*}
has two flaps separated by one tile.

\noindent%
This assemblage of the sequence (\ref{eq:seqf2a6e}, right) has two flaps separated by two tiles:
\begin{equation*}
\sla E_1 E N W_2 N S_1 S_2 W .
\end{equation*}

\begin{figure}\begin{center}
\includegraphics[width=3cm]{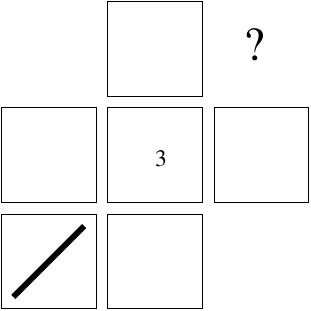}~~~
\includegraphics[width=3cm]{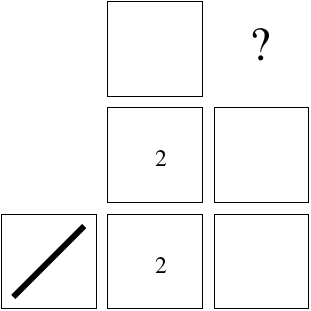}~~~
\includegraphics[width=3cm]{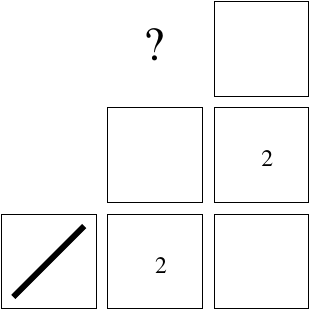}
\caption{\small{Schematic structure of the three unclassified sequences with two flaps.}}\label{fig:unclassified1}
\end{center}\end{figure}

\noindent%
Finally there is one unclassified assemblage with two tiles separated by three tiles:
\begin{equation}
\sla E_1 E N_1 N S_2 W S_2 W .
\end{equation}
A schematic visualization of these assemblages is depicted in Figure \ref{fig:unclassified1}

%
\begin{figure}\begin{center}
\includegraphics[width=3cm]{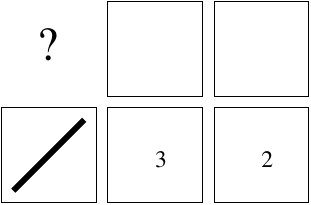}~~~
\includegraphics[width=3cm]{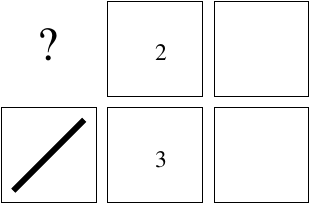}
\caption{\small{Schematic structure of the two unclassified sequences with three flaps.}}\label{fig:unclassified2}
\end{center}\end{figure}

\subsubsection{Sequences with three \textit{flaps}}

There are 2 of them:
\begin{eqnarray*}
\sla E_3 E_2 W_2 E_1 N W S_1 W
\\
\sla E_3 E N W_1 S_1 N_2 S_2 W .
\end{eqnarray*}
The first one would conceivably be constructible starting from 
$\sla E_3 E_2 W_1 E_1 N W S_2 W$ by exchanging the position of two
tiles, which is possible only with a very large amount of stretching on
the wires and deformation of the tiles, not available on the real puzzle.
A schematic visualization of these assemblages is depicted in Figure \ref{fig:unclassified2}

\section{The software codes}\label{sec:code}
The software code is contained in the Subversion (svn) repository \cite{mine:svn},
together with a povray module that can be used to produce syntetized images of
3D configurations and scripts to produce a printout with customized decorations
for the puzzle.
The code can also be downloaded from \cite{mine:web} and should
work on any computer with a \texttt{C} compiler.

There are actually two distinct software codes, described in
Section \ref{sec:code3D} (octominoid shapes) and Section \ref{sec:code2D}
(planar face-up configurations).

\subsection{3D octominoid shapes}\label{sec:code3D}

The name of the executable is \texttt{rubiksmagic}.
If run without arguments, it will search for all canonical representatives
of the set of octominoid configurations that satisfy the local constraints (feasible
configurations), using the
\textit{magic code} described in Section \ref{sec:octominoid}.

This is part of its output:

\begin{Verbatim}[fontsize=\small,frame=single,xleftmargin=10mm,commandchars=\\\{\}]
$ \textit{./rubiksmagic}
RRRURRRU box=0x2x4+8 polyominoid=011-013-015-017-031-033-035-037s
         f=0 delta=0 linking=0 symcount=4 typeinv=yes
RRRUmRRRUm box=1x1x4+8 polyominoid=011-013-015-017-101-103-105-107s
         f=0 delta=0 linking=0 symcount=4 typeinv=no
[...]
RmUvRvLmUvUmUvLv box=2x2x2+0 polyominoid=112-...-312
         f=4 delta=4 linking=0 symcount=1 typeinv=no
RmUvUmLvRmUvUmLv box=2x2x2+0 polyominoid=112-...-323s
         f=4 delta=0 linking=0 symcount=4 typeinv=yes
Found 1291 sequences
$
\end{Verbatim}

The output actually consists of a single (long) line for each configuration, here wrapped in two lines for
convenience, and contains the following information.

\begin{description}
 \item[box]
  The bounding box of the octominoid shape with syntax \texttt{$x$x$y$x$z$+$t$} where $x$, $y$, $z$ are the
  dimensions of the smallest box that contains the shape (rotated such that $x \leq y \leq z$.
  The special boxes \texttt{0x2x4} and \texttt{0x3x3} correspond to the two basic flat octominoid configurations,
  apart from these, we have exactly eight distinct boxes:
  \texttt{1x1x2} (colored purple in \cite{basteleien:web}),
  \texttt{1x1x3} (green),
  \texttt{1x1x4} (blue),
  \texttt{1x2x2} (red),
  \texttt{1x2x3} (grey),
  \texttt{1x3x3} (dark grey),
  \texttt{2x2x2} (light blue),
  \texttt{2x2x3} (yellow).
  The number following the \texttt{+} sign $0 \leq t \leq 8$ denotes the number of tiles of the configuration that lie at the
  boundary of the bounding box;
 \item[polyominoid]
  Describes the actual 3D shape, there is a sequence of eight groups of triplets of digits.
  Each group describes the 3D position of a tile as follows: observe that in each group exactly one of the three
  digits is even, after dividing all digits by two we obtain the 3D coordinates of the center of the tile.  This
  uniquely determines the orientation of the tile by observing that the only integral coordinate indicates the direction
  of the normal vector to the tile;
 \item[f]
  The number of flaps in the configuration;
 \item[delta]
  This is the computed value of the metric invariant introduced in Section \ref{sec:invariantemetrico};
 \item[linking]
  This is the computed linking number $L$ (topological invariant) introduced in Section \ref{sec:invariantetopologico};
 \item[symcount]
  The order of the group of symmetries.
  For example, a configuration that is mirror symmetric will have \texttt{symcount} at least two
  (Figure \ref{fig:shapes}, center, is an example).  The possible values are
  $1$, $2$, $4$, $8$.
  The configuration of Figure \ref{fig:shapes} (right) is completely unsymmetric, hence \texttt{symcount} is $1$.
  There are $4$ configurations with the maximal value $8$ for symcount, but only one of these: \texttt{RRmRRmRRmRRm}
  has both vanishing invariants and has the shape of the lateral surface of a square prism of side $2$ and height $1$.
  Because of the marked diagonals of the tiles, this configuration is \textbf{not} symmetric with respect to reflection
  about an intermediate plane parallel to the base of the prism, thus reducing the number of symmetries from $16$
  (symmetries of a square prism) to $8$;
 \item[typeinv]
  Usually, exchanging the type of all the tiles from $\sla$ to $\bsla$ and viceversa produces a non-equivalent configuration
  (\texttt{typeinv=no}).
  However, $141$ of the $1291$ feasible configurations produce an equivalent configuration upon such exchange
  (\texttt{typeinv=yes}), $50$ of them have vanishing invariants.
  One of these is shown in Figure \ref{fig:shapes} (left).
\end{description}

It is clear from the description above that the software code is able to compute both invariants (in contrast with the
code for the planar face-up configurations described in the next subsection, for which the code is currently not able
to compute the linking number).

Using option \texttt{-M} (\texttt{./rubiksmagic -M}) has the effect of filtering out all configurations with nonvanishing
metric invariant ($\Delta \neq 0$).
Similarly option \texttt{-T} filters out configurations with nonvanishing topological invariant ($L \neq 0$).

Option \texttt{-w} allows to display the warp code described by Verhoeff in \cite{Ver:87};
a few other options are described by \texttt{./rubiksmagic --help}.

The \texttt{rubiksmagic} command accepts an argument, in the form of a magic code, in which case the computation is
limited to the corresponding configuration:
the code first computes the canonical (equivalent) magic code and displays all the above informations for that
configuration.

Finally, the code is not limited to the case of the puzzle with eight tiles;
option \texttt{-n $n$} can be use to compute with the puzzle with $n$ tiles ($n$ must be even).

\begin{figure}\centering
\includegraphics[height=45mm]{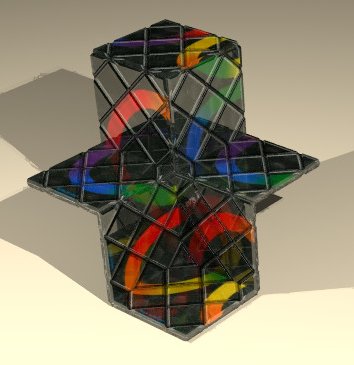}
\includegraphics[height=45mm]{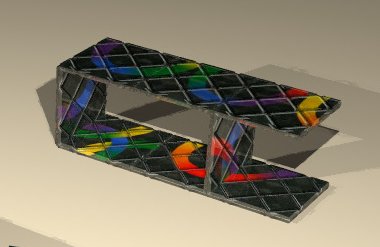}
\caption{\small{The configuration on the left has magic code \texttt{RmRvUvDmLvLmDmUv}.
The configuration on the right has magic code \texttt{RRRmRmRRLmLm}.
These computer generated images were obtained from the magic codes of the configurations using the
PovRay raytracing program and the include file provided in \cite{mine:svn}.}}
\label{fig:nc_shapes}
\end{figure}

As an example, the configuration having magic code \texttt{RmRvUvDmLvLmDmUv} (Figure \ref{fig:nc_shapes}, left)
is analyzed as

\begin{Verbatim}[fontsize=\small,frame=single,xleftmargin=10mm,commandchars=\\\{\}]
$ \textit{./rubiksmagic -c RmRvUvDmLvLmDmUv}
RmRvUvDmLvLmDmUv box=2x2x2+2 polyominoid=011-112-121-...-433s
                 f=0 delta=-4 linking=0 symcount=2 typeinv=yes
\end{Verbatim}
where option \texttt{-c} inhibits the canonization of the magic code.
Since $\Delta \neq 0$ such configuration cannot be constructed.

The configuration with magic code \texttt{RRRmRmRRLmLm} (Figure \ref{fig:nc_shapes}, right)
is analyzed as

\begin{Verbatim}[fontsize=\small,frame=single,xleftmargin=10mm,commandchars=\\\{\}]
$ \textit{./rubiksmagic -c RRRmRmRRLmLm}
RRRmRmRRLmLm box=1x1x3+7 polyominoid=011-013-015-110-114-211-213-215s
             f=2 delta=0 linking=1 symcount=2 typeinv=yes
\end{Verbatim}
and is not constructible since the topological invariant $L = 1$ does not vanish.

\subsection{Planar face-up shapes}\label{sec:code2D}

The name of the executable is \texttt{rubiksmagic2d}.
If run without arguments, it will search for all canonical representatives
of the set $\SSS$ of equivalent classes of sequences.

this is part of its output:

\begin{Verbatim}[fontsize=\small,frame=single,xleftmargin=10mm,commandchars=\\\{\}]
$ \textit{./rubiksmagic2D}
EEEEWWWW f=2 area=5 Dc=0 symcount=8 assemblages=1 deltaiszero=1
EEENWWWS f=0 area=8 Dc=0 symcount=4 assemblages=1 deltaiszero=1
EEENWWSW f=1 area=7 Dc=-2 symcount=1 assemblages=2 deltaiszero=1
[...]
ENSWENSW f=4 area=3 Dc=0 symcount=8 assemblages=0 deltaiszero=0
EWEWEWEW f=8 area=2 Dc=0 symcount=32 assemblages=0 deltaiszero=0
Found 71 sequences
$
\end{Verbatim}

It searches for all admissible sequences that are the canonical representative
of their equivalent class in $\SSS$ (it finds $71$ equivalent classes), for
each one it prints the sequence followed by some information (to be explained shortly).

The software allows for puzzles with a different number of tiles, for example for the
large version with $12$ tiles of the puzzle it finds $4855$ equivalence classes,
with a command like

\begin{Verbatim}[fontsize=\small,frame=single,xleftmargin=10mm,commandchars=\\\{\}]
$ \textit{./rubiksmagic2D -n 12}
EEEEEEWWWWWW f=2 area=7 Dc=0 symcount=8 assemblages=1 deltaiszero=1
EEEEENWWWWWS f=0 area=12 Dc=0 symcount=4 assemblages=1 deltaiszero=1
EEEEENWWWWSW f=1 area=11 Dc=-2 symcount=1 assemblages=2 deltaiszero=1
[...]
ENSNSWENSNSW f=8 area=3 Dc=0 symcount=4 assemblages=0 deltaiszero=0
ENSNSWENSWEW f=8 area=3 Dc=0 symcount=2 assemblages=0 deltaiszero=0
ENSWENSWENSW f=6 area=3 Dc=0 symcount=12 assemblages=0 deltaiszero=0
EWEWEWEWEWEW f=12 area=2 Dc=0 symcount=48 assemblages=0 deltaiszero=0
Found 4855 sequences
$
\end{Verbatim}
however the computational complexity grows exponentially with the number of tiles.

Another use of the code allows to ask for specific properties of a given sequence,
we illustrate this with an example:

\begin{Verbatim}[fontsize=\small,frame=single,xleftmargin=10mm,commandchars=\\\{\}]
$ \textit{./rubiksmagic2D -c EEWENWSW}
EEWENWSW f=3 area=5 Dc=-2 symcount=1 assemblages=6 deltaiszero=2
 Assemblage with delta = 0: sla E3 E2 W2 E1 N1 W1 S1 W1
 Assemblage with delta = 0: sla E3 E2 W1 E1 N1 W1 S2 W1
$
\end{Verbatim}

The first line of output displays some information about the sequence given in the
command line, specifically we find
\begin{itemize}
\item the sequence itself;
\item the number of flaps ($3$ in this case);
\item the area of the plane covered ($5$);
\item the computed contribution $\Delta_c$ coming from the curving tiles;
\item the cardinality of the group of symmetries of the sequence, this particular
sequence does not have any symmetry;
\item the number of admissible assemblages of the sequence, counting only those that
start with $T_0$ of type $\sla$ and identifying assemblages that are equivalent under
transformations in the group of symmetries of the sequence;
\item the number of admissible assemblages with vanishing metric invariant
($\Delta = 0$), we have two in this case.
\end{itemize}

Then we have one line for each of the possible assemblages with $\Delta = 0$ with
a printout of each assemblage, the numbers after each cardinal direction tells the
level of the tile reached with that direction.
It will be $1$ for tiles that are not superposed with other tiles, otherwise it is
an integer between $1$ and the number of superposed tiles.

The option `\texttt{-c}' on the command line can be omitted in which case the
software computes the canonical representative of the given sequence and prints
all the informations for both the original sequence and the canonical one.
Note that the sign of the invariants is sensitive to equivalence transformations.




\begin{thebibliography}{99}
%

\bibitem{wikipedia:web} Rubik's Magic - Wikipedia, 
\url{https://en.wikipedia.org/wiki/Rubik's_Magic},
retrieved Jan 15, 2014.

\bibitem{basteleien:web} J. K\"oller, Rubik's Magic,
\url{http://www.mathematische-basteleien.de/magics.htm},
retrieved Jan 15, 2014.

%
\bibitem{Nou:86} J.G. Nourse,
Simple Solutions to Rubik's Magic,
New York, 1986.

\bibitem{mine:web} M. Paolini, Rubik's Magic,
\url{http://rubiksmagic.dmf.unicatt.it/},
retrieved Jan 15, 2014.

\bibitem{mine:svn} M. Paolini, rubiksmagic project,
(2015), Subversion repository,
\url{https://svn.dmf.unicatt.it/svn/projects/rubiksmagic/trunk}.

%

%
\bibitem{jaap:web} Jaap Scherphuis, Rubik's Magic Main Page,
\url{http://www.jaapsch.net/puzzles/magic.htm},
retrieved Jan 15, 2014.

%
\bibitem{Ver:87} Tom Verhoeff,
{\it Magic and Is Nho Magic},
Cubism For Fun {\bf 15} (1987), 24--31.

%
\end{thebibliography}
\end{document}